\documentclass[10pt,a4paper]{amsart}

\usepackage{amsmath,amssymb}

\input xy
\xyoption{all}

 at10pt

\def\ext{\mathop{\rm Ext}\nolimits}
\def\tor{\mathop{\rm Tor}\nolimits}

\def\spec{\mathop{\rm Spec}}
\def\maxspec{\mathop{\rm maxSpec}}
\def\proj{\mathop{\rm Proj}}
\def\codim{\mathop{\rm codim}}
\def\depth{\mathop{\rm depth}}
\def\reg{\mathop{\rm reg}\nolimits}
\def\nreg{\mathop{\rm Sing}}

\def\a{{\alpha}}
\def\b{{\beta}}
\def\g{{\gamma}}

\def\fin{\mathop{\rm end}}
\def\pd{\mathop{\rm pdim}}
\def\chdim{\mathop{\rm cd}}

\def\coker{\mathop{\rm coker}}

\def\codim{\mathop{\rm codim}}
\def\indeg{\mathop{\rm indeg}}

\def\depth{\mathop{\rm depth}}

\def\supp{\mathop{\rm Supp}}
\def\sym{\mathop{\rm Sym}\nolimits}

\def\D{{\mathcal D}}
\def\C{{\mathcal C}}
\def\F{{\mathcal F}}

\def\Z{{\mathcal Z}}

\def\sh{{\mathcal S}}

\def\om{{\omega}}
\def\Om{{\Omega}}

\def\a{{\alpha}}
\def\e{{\varepsilon}}
\def\b{{\beta}}
\def\g{{\gamma}}
\def\t{{\theta}}
\def\d{{\delta}}

\def\ip{{\mathfrak p}}

\def\idd{{\mathfrak d}}
\def\im{{\mathfrak m}}

\def\ra{{\rightarrow}}
\def\lra{{\longrightarrow}}
\def\fini{{$\quad\quad\Box$}}

\newtheorem{thm}{Theorem}[section]
\newtheorem{lemma}[thm]{Lemma}
\newtheorem{cor}[thm]{Corollary}
\newtheorem{pro}[thm]{Proposition}
\newtheorem{example}[thm]{Example}
\newtheorem{remark}[thm]{Remark}
\newtheorem{ddef}[thm]{Definition}

\newcommand{\bd}{\begin{ddef}}
\newcommand{\ed}{\end{ddef}}
\newcommand{\bt}{\begin{thm}}
\newcommand{\et}{\end{thm}}
\newcommand{\bl}{\begin{lemma}}
\newcommand{\el}{\end{lemma}}
\newcommand{\bco}{\begin{cor}}
\newcommand{\eco}{\end{cor}}
\newcommand{\bp}{\begin{pro}}
\newcommand{\ep}{\end{pro}}
\newcommand{\bex}{\begin{example}}
\newcommand{\eex}{\end{example}}
\newcommand{\brm}{\begin{remark}}
\newcommand{\erm}{\end{remark}}
\newcommand{\bconj}{\begin{conj}}
\newcommand{\econj}{\end{conj}}

\newcommand{\bib}{\bibitem}

\newcommand{\beqn}{\begin{eqnarray*}}
\newcommand{\eeqn}{\end{eqnarray*}}
\newcommand{\beq}{\begin{eqnarray}}
\newcommand{\eeq}{\end{eqnarray}}
\newcommand{\been}{\begin{enumerate}}
\newcommand{\eeen}{\end{enumerate}}

\begin{document}

\author[M. Chardin]{Marc Chardin}
\title[Castelnuovo-Mumford regularity \& some functors]
{On the behavior of Castelnuovo-Mumford regularity with respect to some functors}

\address{M. Chardin, Institut Mathematiques de Jussieu,
175 rue du Chevaleret,
F-75013, Paris, France.}
\email{chardin@math.jussieu.fr}




\maketitle

\begin{abstract}
We investigate the behavior of Castelnuovo-Mumford regularity with respect to 
some classical functors : Tor, the Frobenius functor in positive characteristic, taking a power
or a product (on ideals). These generalizes and refines previous results on these issues by several authors. 
As an application we provide results on the regularity of an intersection of subschemes of a 
projective scheme, under appropriate geometric hypotheses.  
Results on the rigidity of multiple Tor modules and on the characterization of their vanishing are given,
motivated by geometric applications.
\end{abstract}

\section*{Introduction}

In this article we investigate the behavior of Castelnuovo-Mumford regularity with respect to 
some classical functors. One of our first motivations was to provide estimates on
the regularity in a geometric context. This is for instance the content of Theorem \ref{regintersect}
which implies the following result :
\bt
 Let $k$ be a field, $\Z_{1},\ldots ,\Z_{s}$ be closed subschemes of a closed subscheme $\sh\subset {\bf P}^{n}_{k}$. Assume that $\sh$ is irreducible with a singular locus of dimension at most 1. If $\Z :=\Z_{1}\cap \cdots
\cap \Z_{s}\subset \sh$ is a proper intersection of  subschemes of $\sh$ that are Cohen-Macaulay locally
at points of $\Z$, then setting $r'_\sh :=\max\{ \reg (\sh )-1,0\}$, one has
$$
\reg (\Z )\leq \sum_{i=1}^{s}\max\{ \reg (\Z_{i}),r'_\sh \} +\lfloor (\dim \sh -1)/2\rfloor r'_\sh .
$$
In particular, if $\reg (\sh )\leq 1$, then $\reg (\Z )\leq \sum_{i=1}^{s}\reg (\Z_{i})$.
\et

Another motivation was the estimates obtained by Eisenbud, Huneke and Ulrich
on the regularity of Tor modules over a polynomial ring, and its application to estimate the regularity of the
powers of an ideal. Their work was inspired by previous results on the regularity of product of ideals and of tensor product of
modules by Conca and Herzog \cite{CH}, Sidman \cite{Si} and Caviglia \cite{Ca}. They have proved in \cite[3.1]{EHU} an upper bound for the regularity of $\tor_i^R(M,N)$ in terms
of $\reg (M)$ and $\reg (N)$, when $R$ is a polynomial ring and $\tor_1^R(M,N)$ is supported in dimension at most 1. We here extend this result from the class of polynomial rings over a field, to the case of a standard graded
algebra over a Noetherian local ring, in the more general context of multiple Tor modules, provided at least all the modules but one are of finite projective dimension --this assumption cannot be dropped. We moreover show that the estimate is 
sharp, in the sense that at least for one $i$ the estimate on the regularity of the $i$-th Tor is an equality. In the case 
of two modules it gives :
\bt\label{regtorint}
Let $S$ be a standard graded ring over a Noetherian local ring $(S_0,\im_0 )$ and $M,N$
be finitely  generated graded $S$-modules.  If $M$ or $N$ has finite projective dimension  and 
$\tor_{i}^{S}(M,N)\otimes_{S_0}S_0/\im_0$ is supported in dimension at most one for $i\geq 1$,  then 
$$
\max_{i}\{ \reg (\tor_{i}^{S}(M,N))-i\} =\reg (M)+\reg (N)-\reg (S).
$$
\et

An independant proof of this result in the case where $N=S_0$ is a field was given by R\"omer in \cite{Roe} in the  
more general setting of positively graded algebras over a field. 
The hypothesis on the dimension of the support of all positive Tor modules (in place of only the first,
in the polynomial case) is needed since Tor is not rigid when $S$ is singular. It however can be 
dropped if $\dim \nreg (S)\leq 1$.

We review some properties of multiple Tor module in the first section, and prove several facts that does
not seem to have been revealed before, in particular the rigidity of multiple Tor modules over a regular ring containing a field (I don't know if this last assumption can be dropped) and the following result that is useful for geometric applications :
 \bt  Let $R$ be a regular local ring containing a field,
$M_{1},\ldots ,M_{s}$ be finitely generated $R$-modules. The following
are equivalent,\smallskip

{\rm (i)}  $\tor_{1}^{R}(M_{1},\ldots ,M_{s})=0$ and 
$M_{1}\otimes_{R}\cdots\otimes_{R}M_{s}$ is Cohen-Macaulay,\smallskip

{\rm (ii)} the codimension of $M_{1}\otimes_{R}\cdots\otimes_{R}M_{s}$ is 
the sum of the projective dimensions of the $M_{i}$'s,\smallskip

{\rm (iii)} the intersection of the $M_{i}$'s is proper
and every $M_{i}$ is Cohen-Macaulay.\smallskip
\et

Over a field $k$ of  characteristic $p>0$, an interesting question was raised by Katzman in \cite{Ka}, 
 motivated by the open question of commutation of localization with tight closure. 
Denote by 
$\F$ the Frobenius functor and by $\F^e$ its iteration $e$ times. Let $M$ be a finitely generated graded module
over a standard graded algebra $S$ over $k$.  Does there exists $C$ such that $\reg (\F^e M)\leq C p^e$ ?

Katzman conjectured that the answer is always yes, and proved it in some cases (monomial ideals in a toric ring). 
We give a partial answer to this question, with hypotheses of different nature. Our key estimate is in terms of 
graded Betti numbers of $M$ over $S$. The following result follows from it :

\bt
Let $S$ be a standard graded ring over a field 
of characteristic $p>0$ and $M$ a finitely generated graded $S$-module. If 
$\dim (\nreg (S)\cap \supp (M))\leq 1$, and $S$ is not regular, then
$$
 \reg  (\F^{e}M)\leq \reg (S)+p^{e}\left( \reg (M)+\left\lfloor {{\dim S}\over{2}}\right\rfloor (\reg (S)-1)+\dim S\right) .
 $$
\et

(If $S$ is regular, $\F$ is exact and $ \reg  (\F^{e}M)\leq p^{e}\left( \reg (M)+\dim S\right)$.)\medskip

When $I$ is an homogeneous ideal in a polynomial ring $R$ over a field such that $\dim R/I$ is
of dimension at most 1, it has been proved by Chandler and Geramita, Gimigliano and Pitteloud in \cite{Chan} and \cite{GGP} that 
$\reg (I^j)\leq j\reg (I)$ for any $j$. This has been refined by Eisenbud, Huneke and Ulrich who 
proved in \cite{EHU} that in this situation $\reg (I^j)\leq \reg (I)+(j-1)(e-1)$ if $I$ is generated in degrees at
most $e-1$ and related in degrees at most $e$. We first show a second refinement of the initial estimate, 
which holds for an homogeneous ideal $I$ of a Noetherian standard graded ring $S$ : 

\bt\label{regpow1int}
Let $I$ be an homogeneous $S$-ideal  such that $\dim (S/I)\otimes_{S_0}S_0/\im_0\leq 1$ for any maximal ideal $\im_0\in \spec (S_0)$. Then, for any $m\geq 0$,
$$
\reg (S/I^{m+1})\leq \max\{ a_{0}(S/I)+b_{0}^{S}(I),(a_{1}(S/I)+1)+\reg_{1}^{S}(I)\} +(m-1)b_{0}^{S}(I).
$$
\et
Notice that $a_{1}(S/I)+1$ is the regularity of $S/I^{sat}$, where $I^{sat}$ is the saturation of $I$ 
with respect to the positive part of $S$.

Second we extend this result to dimension 2, under the assumption that $I$ is generically a complete 
intersection in $S$ and $S_0$ is a field. Recall that Terai and Sturmfels provided examples such that $\reg (I^2)>2\reg (I)$, so
that one cannot hope for straightforward extension of this result. In the particular case of the polynomial ring,
our result implies the following 

\bt  Let $I$ be an homogeneous ideal of a polynomial ring 
  $R$ over a field, generated in degrees at most $d$, such that $\dim (R/I)=2$. Assume that 
  $I_{\ip}\subset R_{\ip}$ is a complete intersection  for every prime $\ip\supseteq I$ such that $\dim R/\ip =2$. 
  
  Then $\reg (I^{2})\leq \max \{ 2\reg (I),\reg (I^{sat})+2d-2\}$ and, for $j\geq 3$,
$$
\reg (I^{j})\leq  \max \{ 3\reg (I)+(j-3)d,\reg (I^{sat})+jd-2\} .
$$
\et

One of the central tools in our study is the use of non acyclic complexes to control the Castelnuovo-Mumford 
regularity of a module. This idea was introduced by Gruson, Lazarsfeld and Peskine to bound the regularity
of reduced curves in \cite[1.6]{GLP}. It has been exploited by many authors since, particularly in the recent years.

We provide in section 3 several lemmas that exploits  this technique in several contexts, and with different level of
refinements. Their proofs all rely on the analysis of a spectral sequence with different
types of degeneration. 

The second section gives standard estimates on graded Betti numbers. They all rely on techniques that
have been used in closely related context. We have provided complete (however short) proofs for lack
of a reference that directly applies to our situation. 

\medskip

The sections are organized as follows :

\medskip

\S 1. Some results on Tor modules

\S 2. Standard estimates on graded Betti numbers

\S 3. Non acyclic complexes and Castelnuovo-Mumford regularity

\S 4. Regularity of Frobenius powers

\S 5. Regularity of Tor modules

\S 6. On the regularity of ordinary powers

\bigskip
\section{Some results on Tor modules}
\medskip

\bd\label{defmtor} Let  $A$ be a ring, $M_{1},\ldots
,M_{s}$ be $A$-modules and $F_{\bullet}$ be  the tensor products 
over $A$ of the canonical free resolutions  of $M_{1},\ldots ,M_{s}$. Then 
$$
\tor_{i}^{A}(M_{1},\ldots ,M_{s}):=H_{i}(F_{\bullet}).
$$
\ed

Notice that if $F_{\bullet}'$ is the product of any choice of flat resolutions of  
$M_{1},\ldots ,M_{s}$, then $\tor_{i}^{A}(M_{1},\ldots ,M_{s})\simeq 
H_{i}(F_{\bullet}')$.
\medskip
\bp\label{A.2} Let $A$ be a ring, $M$ be an $A$-module and $F_{\bullet}$ be a 
complex of flat $A$-modules. If $L_{\bullet}$ is a resolution of $M$, there is a 
natural isomorphism, 
$$
H_{i}(F_{\bullet}\otimes_{A}L_{\bullet})\simeq H_{i}(F_{\bullet}\otimes_{A}M).
$$
\ep

{\it Proof.} Consider the spectral sequence with
$E^{1}_{ij}=H_{j}(F_{i}\otimes_{A}L_{\bullet})=F_{i}\otimes_{A}H_{j}(L_{\bullet})$
and $E^{2}_{ij}=H_{i}(E^{1}_{\bullet j})$ that 
abouts to $H_{i}(F_{\bullet}\otimes_{A}L_{\bullet})$ and note that by
hypothesis $H_{j}(L_{\bullet})=0$ for $j\not= 0$ and
$H_{0}(L_{\bullet})=M$.\fini\medskip

\bco\label{A.3} Let $A$ be a ring, $M,M_{1},\ldots
,M_{s}$ be $A$-modules and $F_{\bullet}$ be  the tensor products 
over $A$ of free resolutions of $M_{1},\ldots ,M_{s}$. Then,
$$
\tor_{i}^{A}(M_{1},\ldots ,M_{s},M)\simeq H_{i}(F_{\bullet}\otimes_{A}M).
$$
\eco

If $(R,\im )$ is local regular of dimension $n$, this result implies that 
$\tor_{i}^{R}(M_{1},\ldots ,M_{s})=0$ for $i>n(s-1)$. Also notice in this case 
that if $M_{1}=\cdots =M_{s}=R/\im$, $\tor_{n(s-1)}^{R}(M_{1},\ldots ,M_{s})\simeq
R/\im \not= 0$.
\medskip

\bp\label{A.4} Let $A$ be a ring and $M_{1},\ldots
,M_{s},N_{1},\ldots ,N_{t}$ be $A$-modules. There exists a
spectral sequence,
$$
E^{2}_{ij}=\tor_{i}^{A}(M_{1},\ldots ,M_{s},\tor_{j}^{A}(N_{1},\ldots
,N_{t}))\ 
\Rightarrow \ \tor_{i+j}^{A}(M_{1},\ldots ,M_{s},N_{1},\ldots
,N_{t}).
$$
\ep

{\it Proof.} Let $F_{\bullet}$ and $L_{\bullet}$ be  the tensor products 
over $A$ of free resolutions of $M_{1},\ldots ,M_{s}$ and $N_{1},\ldots ,N_{t}$, 
respectively. 

The double complex $C$ with $C_{ij}:=F_{i}\otimes_{A}L_{j}$ gives 
rise to two spectral sequences both abouting to $H_{\bullet}({\rm Tot}(C))\simeq \tor_{\bullet}^{A}(M_{1},\ldots ,M_{s}
,N_{1},\ldots,N_{t})$. One of them have as second terms
$E^{2}_{ij}=H_{i}(F_{\bullet}\otimes_{A}H_{j}(L_{\bullet}))$
and the result follows from Corollary \ref{A.3}.\fini\medskip

Other spectral sequences for multiple Tor modules appears in  \cite[\S 6]{EGA3}.

\bp\label{A.5} If $B\lra A$ is a flat map and 
$M_{1},\ldots,M_{s}$ are $B$-modules, then
$$
\tor_{i}^{A}(M_{1}\otimes_{B}A,\ldots ,M_{s}\otimes_{B}A)\simeq
\tor_{i}^{B}(M_{1},\ldots ,M_{s})\otimes_{B}A
$$
for every $i$.
\ep

{\it Proof.} This is a special case of \cite[6.9.2]{EGA3}. \fini\medskip

\brm\label{A.6} If $I_{1},\ldots ,I_{s}$ are ideals in a ring $A$, then
$\tor_{1}^{A}(A/I_{1},\ldots ,A/I_{s})$ admits the following
description that generalizes the isomorphism
$\tor_{1}^{A}(A/I,A/J)\simeq (I\cap J)/IJ$. 

Let $M$ be the submodule of $I_{1}\oplus \cdots \oplus I_{s}$ of
tuples $x=(x_{1},\ldots ,x_{s})$ such that $x_{1}+\cdots +x_{s}=0$. Then
$M$ contains the submodule $P$ generated by the tuples $x$ such that
$x_{\ell}=0$ except for two indices $i$ and $j$ and $x_{i}=-x_{j}\in I_{i}I_{j}$.
Then,
$$
\tor_{1}^{A}(A/I_{1},\ldots ,A/I_{s})\simeq M/P.
$$
\erm

\bt\label{rigtor} Let $R$ be a regular local ring containing a field
and $M_{1},\ldots ,M_{s}$ be finitely generated
  $R$-modules. Then,\smallskip

{\rm (i)} $\tor_{i}^{R}(M_{1},\ldots ,M_{s})=0$ implies 
  $\tor_{j}^{R}(M_{1},\ldots ,M_{s})=0$  for all $j\geq i$.\smallskip

{\rm (ii)} Let $j:=\max\{ i\ |\ \tor_{i}^{R}(M_{1},\ldots ,M_{s})\not= 0\}$. Then
$$
\pd M_{1}+\cdots +\pd M_{s}=\dim R+j-\varepsilon ,
$$
with $0\leq \varepsilon \leq \dim \tor_{j}^{R}(M_{1},\ldots ,M_{s})$. Let
$\varepsilon_{0}:=\min_{i}\{ \depth  \tor_{j-i}^{R}(M_{1},\ldots ,M_{s})
+i\}$, then $\varepsilon \geq  \varepsilon_{0}$ and  equality holds if $\varepsilon_{0}=\depth  \tor_{j}^{R}
(M_{1},\ldots ,M_{s})$.
\et

{\it Proof.} First we complete $R$ for the $\im$-adic filtration (where $\im$ is the 
maximal ideal of $R$), and use Cohen 
structure theorem to reduce to the case of a power series ring over $k:=R/\im$. Notice
that $R\lra \hat R$ is flat and as the $M_{i}$'s are finite $\hat M_{i}=M_{i}
\otimes_{R}\hat R$, $\dim_{R}M_{i}=\dim_{\hat R}\hat M_{i}$ and  $\depth_{R} 
M_{i}=\depth_{\hat R}\hat M_{i}$ for every $i$. 

We let $n:=\dim R$, and may assume that $R=k[[X_{1},\ldots ,X_{n}]]$.

Consider $S$ the completed tensor product over $k$ of $s$ copies of $R$. Notice
that  $S\simeq k[[X_{ij}\ |\ 1\leq i\leq n,\ 1\leq j\leq s ]]$, if
we set $X_{ij}:=1\otimes \cdots \otimes 1 \otimes
X_{i}\otimes 1\cdots \otimes 1$ where $X_{i}$ is at the $j$-th
position. The ring $S$ is local regular with maximal ideal $\im_{S}=(X_{ij})=
\idd_{S}+\im /\idd_{S}$,  where $\idd_{S}$ defined by the exact sequence
$$
0\lra \idd_{S}\lra S\buildrel{\rm mult}\over\lra R\lra 0.
$$
Notice that $\idd_{S}$ is generated by a regular sequence, for instance 
$f:=(X_{ij}-X_{i(j+1)}\ |\ 1\leq i\leq n,\
1\leq j<s )$ is a minimal generating $n(s-1)$-tuple. 
Considering $R$ as an $S$-module via the above exact sequence, for any tuple
$L_{1},\ldots ,L_{t}$ of free $R$-modules, there is a natural isomorphism
$$
L_{1}\otimes_{R}\cdots \otimes_{R}L_{t}\buildrel{\sim}\over\lra 
(L_{1}\hat\otimes_{k}\cdots \hat\otimes_{k}L_{t})\otimes_{S}R.
$$

We put $P:=M_{1}\hat\otimes_{k}\cdots \hat\otimes_{k}M_{s}$ choose
a minimal free $R$-resolution $F^{(i)}_{\bullet}$  of $M_{i}$, set 
$F^{k}_{\bullet}:=F^{(1)}_{\bullet}\hat\otimes_{k}\cdots
\hat\otimes_{k}F^{(s)}_{\bullet}$ and $F^{R}_{\bullet}:=F^{(1)}_{\bullet}
\otimes_{R}\cdots\otimes_{R}F^{(s)}_{\bullet}$. 

Notice that $F^{k}_{\bullet}$ is a minimal free $S$-resolution of $P$ and
$K_{\bullet}(f;S)$ is a minimal free $S$-resolution of $R=S/\idd_{S}$. 
We therefore have isomorphismes, 
$$
\tor_{i}^{R}(M_{1},\ldots ,M_{s})\simeq H_{i}(F^{R}_{\bullet})
     \simeq H_{i}(F^{k}_{\bullet}\otimes_{S}R)
     \simeq \tor_{i}^{S}(P,R)
     \simeq H_{i}(K_{\bullet}(f;P)),
$$
and (i) follows by \cite[1.6.31]{BH} as $(S,\im_{S})$ is local Noetherian and 
$\idd_{S}\subset \im_{S}$.

For (ii) let $T^{\bullet}$ be the total complex associated to the 
double complex $\C_{\im_{S}}^{\bullet}K_{\bullet}(f;P)$, with $T^{i}:=\bigoplus_{p-q=i}
\C_{\im_{S}}^{p}K_{q}(f;P)$ and choose the upper index as line index. We will
estimate in two ways the number $\t :=\min \{ i\ |\ H^{i}(T^{\bullet})\not= 0\}$.

We have noticed above that the projective dimension of $P$ (over $S$) is equal to the 
sum of the projective dimensions of the $M_{i}$'s (over $R$).
Equivalently the depth $\d$ of $P$ is equal to the sum of the depths of the $M_{i}$'s.

The homology of $T^{\bullet}$ is the aboutment of two spectral sequences associated
to the horizontal and vertical filtrations of the double complex.

For the vertical filtration, one has $^{v}_{1}E^{p}_{q}=H_{\im_{S}}^{p}K_{q}(f;P)
\simeq H_{\im_{S}}^{p}(P)^{{n(s-1)}\choose{q}}$. Therefore  $^{v}_{1}E^{p}_{q}=0$
for $q>n(s-1)$ and for $p<\depth P=\delta$ which shows that $\t \geq \d -n(s-1)$. 
Also  $^{v}_{2}E^{\d}_{n(s-1)}\simeq {^{v}_{\infty}E^{\d}_{n(s-1)}}$ is 
Matlis dual to $\ext^{ns-\d}_{S}(P,\om_{S})/\idd_{S}$ which is not zero by Nakayama's
lemma. We deduce that $\t = \d -n(s-1)$.

The other spectral sequence has second terms  $^{h}_{2}E^{p}_{q}\simeq 
H^{p}_{\im_{S}}(\tor^{R}_{q}(M_{1},\ldots ,M_{s}))$, and an easy computation 
gives  $\varepsilon_{0}-j=\min\{ p-q\ |\ ^{h}_{2}E^{p}_{q}
\not= 0\}$ so that $\t =\min\{ p-q\ |\ ^{h}_{\infty}E^{p}_{q}
\not= 0\}=:\varepsilon -j\geq \varepsilon_{0}-j$ and it is clear form its definition 
that $ \varepsilon_{0}\geq 0$. Also if $\varepsilon_{0}=\depth \tor^{R}_{j}
(M_{1},\ldots ,M_{s})$ 
we have $^{h}_{2}E^{\varepsilon_{0}}_{j}\simeq {^{h}_{\infty}E^{\varepsilon_{0}}_{j}}$ 
which implies that $\t \leq \varepsilon_{0}-j$ and shows that $\varepsilon =\varepsilon_{0}$ in this case.
It remains to prove that $\varepsilon \leq d:=\dim  \tor^{R}_{j}(M_{1},\ldots ,M_{s})$. 

Notice that ${^{h}_{2}E^{d}_{j}}$ is Matlis dual to 
 $E:=\ext^{ns-d}_{S}( \tor^{R}_{j}(M_{1},\ldots ,M_{s}),\om_{S})$ which is a module
of dimension $d$, and  ${^{h}_{\infty}E^{d}_{j}}$ is Matlis dual to 
a submodule of $E$ that coincides with $E$ in dimension $d-1$, because
${^{h}_{\ell}E^{d-\ell}_{j-\ell+1}}$ for $\ell\geq 2$ is Matlis dual to a module 
supported in dimension at most $d-\ell$. \fini\medskip

\bco\label{A.8}  Let $R$ be a regular local ring containing a field,
$M_{1},\ldots ,M_{s}$ be finitely generated $R$-modules. The following
are equivalent,\smallskip

{\rm (i)}  $\tor_{1}^{R}(M_{1},\ldots ,M_{s})=0$ and 
$M_{1}\otimes_{R}\cdots\otimes_{R}M_{s}$ is Cohen-Macaulay,\smallskip

{\rm (ii)} the codimension of $M_{1}\otimes_{R}\cdots\otimes_{R}M_{s}$ is 
the sum of the projective dimensions of the $M_{i}$'s,\smallskip

{\rm (iii)} the intersection of the $M_{i}$'s is proper
and every $M_{i}$ is Cohen-Macaulay.\smallskip
\eco

{\it Proof.}  Let $M:=M_{1}\otimes_{R}\cdots\otimes_{R}M_{s}$.  Krull's 
intersection theorem gives, 
$$
\codim M\leq  \codim M_{1}+\cdots + \codim M_{s}.
$$

If (i) holds,  Krull's intersection theorem, Theorem \ref{rigtor} (i) and (ii) implies that 
$$
\pd M_{1}+\cdots + \pd M_{s} = n-\depth M
=\codim M \leq  \codim M_{1}+\cdots + \codim M_{s}
$$
which implies (ii) and (iii).

If (ii) holds, Krull's intersection theorem implies that (iii) holds. 

Now assume that (iii) holds. We then have $\codim M=\pd M_{1}+\cdots +\pd M_{s}=n+
j-\varepsilon$, with the notations of  Theorem \ref{rigtor}. Therefore $\varepsilon =
\dim M+j$ which implies that $j=0$ and $\varepsilon =\dim M$ as  $\varepsilon 
\leq\dim M$ by Theorem \ref{rigtor}. This proves (i), and conclude the proof.\fini\medskip

\section{Standard estimates on graded Betti numbers}

We present here some estimates on graded Betti numbers that will be useful in the core of
the article. We don't claim any originality, as all these results already  appeared, at least essentially,
in previous work of several authors (in particular of Avramov).

For a standard graded ring $B$ over a field $k$ and a finitely generated graded $B$-module $M$ set
$$
P^{B}_{M}(t,u):=\sum_{i,j}\dim_{k}\tor_{i}^{B}(M,k)_{j}t^{i}u^{j}.
$$
Recall that if $P=\sum_{i,j}a_{i,j}t^{i}u^{j}$ and $Q=\sum_{i,j}b_{i,j}t^{i}u^{j}$, $P\preccurlyeq Q$ means that 
$a_{i,j}\leq b_{i,j}$ for any $i$ and any $j$.
As $M$ is finitely generated, $P^{B}_{M}(t,u)=\sum_{i}(P_{M}^{B})_{i}(u)t^{i}$, where $(P_{M}^{B})_{i}\in {\bf Z}[u,u^{-1}]$. We set
$b_{i}^{B}(M):=\max\{ j\ \vert\ \tor_{i}^{B}(M,k)_{j}\not= 0\}$ if $\tor_{i}^{B}(M,k)\not= 0$ and $b_{i}^{B}(M):=-\infty$ else,
$\reg_{j}^{B}(M):=\max_{i\leq j}\{ b_{i}^{B}(M)-i\}$ and $\reg^{B}(M):=\sup_{j}\{ \reg_{j}^{B}(M)\}$. With these notations, 

\bp\label{A.9} Let $R$ and $S$ be  standard graded algebras over a field and $R\ra S$ 
a finite degree 0 morphism of graded rings. For any graded $S$-module $M$ of finite type, one has :\smallskip

(i) $P_{M}^{S}(t,u)\preccurlyeq {{P_{M}^{R}(t,u)}\over{1-t(P_{S}^{R}(t,u)-1)}}$.\smallskip

(ii) $P_{M}^{R}(t,u)\preccurlyeq P_{M}^{S}(t,u)P_{S}^{R}(t,u)$.\smallskip


(iii) Set $E_{p}^{R}(S):=\max_{\buildrel{j_{1}+\cdots +j_{\ell}=p}\over{\scriptscriptstyle{j_{1}\geq 2,\ldots ,j_{\ell}\geq 2}}}\{ b_{j_{1}-1}^{R}(S)+\cdots +b_{j_{\ell}-1}^{R}(S)\}$ if $p\geq 2$ and $E_{p}^{R}(S):=0$ else. Then, for any $i$,
$$
\begin{array}{rl}
b_{i}^{S}(M)-i&\leq -i+\max_{p}\{ b_{i-p}^{R}(M)+E_{p}^{R}(S)\} ,\\
&\leq -i+\max_{p}\{ b_{i-p}^{R}(M)+p+\max_{0\leq \ell\leq \lfloor p/2\rfloor }\{ \ell\, (\reg_{p-2\ell+1}^{R}(S)-1) \}\} \\
&\leq  \max_{p}\{ \reg_{i-p}^{R}(M)+ \lfloor p/2\rfloor \max\{ 0,\reg_{p-1}^{R}(S)-1\} \}\\
\end{array}
$$

(iv) For any $i$,
$$
-\reg_{i}^{R}(S)\leq \reg_{i}^{S}(M)-\reg_{i}^{R}(M)\leq  \max \{ 0,\lfloor i/2\rfloor ( \reg_{i-1}^{R}(S)-1)\} .
$$
\ep


{\it Proof.} Point (i)  follows from the Eilenberg-Moore type spectral sequence of Avramov for the DG-module
$N:=\tor^{R}_{\bullet}(M,k)$ over the DG-algebra $A:=\tor^{R}_{\bullet}(S,k)$:
$$
{^{1}E}_{pq}=(S^{k}_{p}(A,N)\otimes_{A}k)_{q}
\Rightarrow \tor^{S}_{p+q}(M,k)
$$
with $S^{k}_{p}(A,N):=A\otimes_{k}\underbrace{\tilde{A}\otimes_{k}\cdots\otimes_{k}\tilde{A}}_{p\ {\rm times}}
\otimes_{k}N$ and $\tilde{A} :=\coker (k\ra A)$. This construction respects the internal gradings of $R$ and $S$ and the 
proof goes along the same lines as in \cite[3.3.2]{Av}.


Claim (ii) is a direct consequence of the change of rings spectral sequence 
$\tor_{p}^{S}(M,\tor_{q}^{R}(S,k))\Rightarrow \tor_{p+q}^{R}(M,k)$.

The first inequality in (iii) is a direct application of the inequality in (i). For the second one,
notice that $b_{j_{t}-1}^{R}(S)\leq \reg_{j_{t}-1}^{R}(S)+j_{t}-1$ and that $\reg_{j_{t}-1}^{R}(S)\leq 
\reg_{p-2\ell -1}^{R}(S)$ as $j_{t}\leq p-2(\ell -1)$. The last inequality is clear from the second one. 

Finally, (ii) implies that $b_{i}^{R}(M)\leq \max_{p+q=i} \{ b_{p}^{S}(M)+b_{q}^{R}(S)\}$ which gives 
the left inequality in (iv), and the right one has been proved in (iii).
\fini\medskip

\bex\label{bettideaux} If $I\subset S$ is a  $S$-ideal and $S=R/J$ where $R$ is a polynomial ring,
then 
$$
\begin{array}{rl}
b_{i}^{S}(I)-i&\leq -i+\max_{p}\{ b_{i-p}^{R}(I/J)+E_{p}^{R}(S)\} ,\\
&\leq -i+\max_{p}\{ b_{i-p}^{R}(I/J)+p+\max_{0\leq \ell\leq \lfloor p/2\rfloor }\{ \ell\, (\reg_{p-2\ell+1}^{R}(S)-1) \}\} \\
&\leq  \max_{p}\{ \reg_{i-p}^{R}(I/J)+ \lfloor p/2\rfloor \max\{ 0,\reg_{p-2}^{R}(J)-2\} \}\\
\end{array}
$$
therefore
$$
-\reg_{i}^{R}(S)\leq \reg_{i}^{S}(I)-\reg_{i}^{R}(I/J)\leq  \max \{ 0,\lfloor i/2\rfloor ( \reg_{i-2}^{R}(J)-2)\} .
$$
Moreover  $\reg_{i}^{R}(I/J)\leq  \max \{ \reg_{i}^{R}(I+J),\reg_{i}^{R}(S)\}$ and $\reg (I/J)=\reg (I+J)$ if $\reg (S/I)\geq \reg (S)$. If $S$ is not a polynomial ring it gives for $i>0$,
$$
\reg_{i}^{S}(S/I)\leq  \max \{ \reg_{i}^{R}(S/I),\reg_{i-1}^{R}(S)-1\} +\lfloor (i-1)/2\rfloor ( \reg_{i-2}^{R}(S)-1).
$$
\eex

\bex\label{A.10} If $R$ is a polynomial ring and $\reg (S)=1$, for instance if $S$ is the coordinate ring of a variety of minimal degree,
then for any graded $S$-module $M$ and any $j$, 
$$
\max_{i\leq j}\{ b_{i}^{R}(M)-i\} -1\leq \max_{i\leq j}\{ b_{i}^{S}(M)-i\} \leq \max_{i\leq j}\{ b_{i}^{R}(M)-i\} .
$$
In particular $b_{i}^{S}(M)\leq \reg (M)+i$ for any $i$. 
\eex

Recall that if $S$ is a standard graded ring over a Noetherian ring and $M$ is a finitely generated graded $S$-module, then $\reg (M):=\max_i \{ a_i (M)+i\} $ with 
$a_i (M):=\max_j \{ H^i_{S_+}(M)\not= 0\}$. 

\bp\label{estbetti}
Let $S$ be a standard graded ring over a field $k$, and $M$ be a finitely generated graded
$S$-module of dimension $d$ . Then,
$$
\reg_{i}^{S}(M)\leq \max_{p\leq d}\{ (a_p (M)+p)+\reg_{i+p}^{S}(k)\} \leq \reg (M)+\reg_{d+i}^{S}(k).
$$
\ep

{\it Proof.} Let $F_{\bullet}$ be a minimal graded free $S$-resolution of $k$ and set $\im:=S_+$. The double complex $\C^{\bullet}_{\im}(F_{\bullet}\otimes_S M)$ gives rise to a spectral sequence
$$
H^p_\im (M)\otimes_S F_q\  \Rightarrow \ \tor_{q-p}^{S}(M,k).
$$
It follows that
$
b_{i}^{S}(M)\leq \max_{p}\{ a_p (M)+b_{p+i}^{S}(k)\} =\max_{p\leq d}\{ (a_p (M)+p)+(b_{p+i}^{S}(k)-p-i)\} +i.
$
\fini

\section{Non acyclic complexes and Castelnuovo-Mumford regularity}

In this section, $R$ is a polynomial ring over a commutative Noethrian ring $R_0$, with its standard grading,   and  $D_{\bullet}$ a graded complex of finitely generated $R$-modules with $D_{i}=0$ for
  $i<0$.

The following lemmas give ways for controlling the local cohomology of
$H_{0}(D_{\bullet})$  in terms of local cohomologies of the $D_{i}$'s under
restrictions on the dimension of the modules $H_{i}(D_{\bullet})$ for $i>0$. 

The invariants $a_j (M)$ are defined as in section 2. The quantities $\d_{p}:= \max_{i\geq 0}\{
a_{p+i}(D_{i})\}$ and $\e_{q}:=\max_{i\geq 0}\{
a_{i}(D_{q+i})\}$ are useful in these estimates.\medskip

Lemmas below all derive from the study of the two spectral sequences
coming from the double complex $\C^{\bullet}_{\im}D_{\bullet}$. For
clarifying later references, we first state particular cases that are often
used in the sequel before giving more refined statements. The proofs of
the six lemmas will be given after they are all stated.\medskip  

Let us first recall a classical fact about the cohomological dimension of a finitely generated graded $R$-module $M$ (denoted by $\chdim M$). A proof of this can be found in \cite[2.3]{Bro} for the local case, to
which the reduction is immediate.  

\bp\label{cohdim}
For $\ip\in \spec (R_0)$ let $k(\ip ):=(R_0)_\ip /\ip (R_0)_\ip$. Then,
 $$
\chdim M:= \max\{ i\ \vert\ H^i_{S_+}(M)\not= 0\} =\max_{\ip\in \spec (k)}\{ \dim (M\otimes_k k(\ip ))\}=\max_{\ip\in \maxspec (k)}\{ \dim (M\otimes_k k(\ip ))\} .
 $$
\ep

We set $H_{i}:=H_{i}(D_{\bullet})$ and introduce the following
conditions on the dimension of the  $H_{i}$'s :

$$
\D_{\ell} (\tau ) :\quad \chdim H_{i}\leq \max\{ \tau ,\tau -\ell +i\} ,\ \forall
i>0.
$$

\bl\label{1.1} $a_{p}(H_{0})\leq \d_{p}$ for $p\leq \tau -1$ if
  $\D_{1}(\tau )$ is satisfied.
  \el
  
\bl\label{1.2}   If $\D_{1} (\tau )$ is
  satisfied, there exists a natural map $\xi :
  H^{\tau -2}_{\im}(H_{0})\lra H^{\tau}_{\im}(H_{1})$ such
  that $\xi_{\mu}$ is onto for $\mu >\d_{\tau +1}$ 
and  $\xi_{\mu}$ is into if $\D_{2} (\tau )$ is
  satisfied and  $\mu >\d_{\tau}$.
  \el
  
  \bl\label{1.3}  If  $\D_{2} (\tau )$ is
  satisfied, there exists a diagram of natural maps,
$$
\xymatrix{
 &
 &0\ar[r]&K\ar[r]\ar^{\t}[d]&H^{\tau -3}_{\im}(H_{0})\ar^{\a}[r]&
H^{\tau -1}_{\im}(H_{1})\\
H^{\tau -4}_{\im}(H_{0})\ar^{\g}[r]&H^{\tau -2}_{\im}(H_{1})\ar^{\b}[r]&
H^{\tau}_{\im}(H_{2})\ar[r]&C\ar[r]&0&\\
}
$$
with $K=\ker (\a )$ and $C=\coker (\b )$, such that,\smallskip

{\rm (a)} $\a_{\mu}$ and  $\t_{\mu}$ are onto for $\mu >\d_{\tau}$,\smallskip

{\rm (b)}  $\t_{\mu}$ is into and ${\rm Im}(\g_{\mu})=\ker (\b_{\mu})$
if $\D_{3} (\tau )$ is   satisfied  and  $\mu >\d_{\tau -1}$.
\el

\bl\label{1.4} Assume that $\chdim H_{i}\leq 1$ for $i\geq 1$. Then 
  $a_{p}(H_{0})\leq \d_{p}$ for $p\geq 0$. Also
  $a_{0}(H_{q})\leq \varepsilon_{q}$ and $a_{1}(H_{q})\leq
  \varepsilon_{q-1}$  for $q>0$. 
  \el
  
 \bl\label{1.5} If $\chdim H_{i}\leq 2$ for $i\geq 2$, then $a_{1}(H_{q})\leq
  \varepsilon_{q-1}$  for $q>0$ and there exists natural maps 
$$
\a^{p} : H^{p}_{\im}(H_{0})\lra H^{p+2}_{\im}(H_{1})\quad \hbox{and}
\quad \phi^{q} : H^{0}_{\im}(H_{q})\lra H^{2}_{\im}(H_{q+1})
$$ 
for $p\geq 0$ and for $q>0$ such that\smallskip 

{\rm (a)} $\a^{p}_{\mu}$ is onto for $\mu >\d_{p+1}$ and is into if 
$\mu >\d_{p}$,\smallskip 

{\rm (b)} $\phi^{q}_{\mu}$ is onto for 
$\mu >\varepsilon_{q-1}$ and is into if $\mu >\varepsilon_{q}$.
\el

\bl\label{1.6}  If $\chdim H_{i}\leq 3$ for $i\geq 3$, then there exist complexes with natural maps
$$
0\ra K_{p}\ra H^{p}_{\im}(H_{0})\buildrel{\a^{p}}\over{\lra} H^{p+2}_{\im}(H_{1})
\buildrel{\b^{p}}\over{\lra} H^{p+4}_{\im}(H_{2})\ra C_{p}\ra 0
$$
for $p\geq -1$ that are exact except possibly in the middle, maps 
$\t^{p}:K_{p}\lra C_{p-1}$ for $p\geq 0$, and for $q>0$ a collection of maps
$$
H^{0}_{\im}(H_{q})\buildrel{\phi^{q}}\over{\lra} H^{2}_{\im}(H_{q+1}),
\ H^{1}_{\im}(H_{q})\buildrel{\psi^{q}}\over{\lra} H^{3}_{\im}(H_{q+1})
\ \hbox{and} \ \ker (\phi^{q})\buildrel{\eta^{q}}\over{\lra} \coker (\psi^{q+1})
$$
(with $\b^{-1}=\psi^{1}$)  such that:\smallskip 

{\rm (a)} $\ker (\b^{p}_{\mu})={\rm im}(\a^{p}_{\mu})$ for  $\mu >\d_{p+1}$,\smallskip

{\rm (b)} $\t^{p}_{\mu}$ is onto for $\mu >\d_{p+1}$ and is into if 
$\mu >\d_{p}$,\smallskip 

{\rm (c)} $\phi^{q}_{\mu}$ is onto for 
$\mu >\varepsilon_{q-1}$,\smallskip

{\rm (d)}  $\psi^{q}_{\mu}$ is into for 
$\mu >\varepsilon_{q-1}$,\smallskip

{\rm (e)}  $\eta^{q}_{\mu}$ is onto for 
$\mu >\varepsilon_{q-1}$ and is into if $\mu >\varepsilon_{q}$.

\el

We now turn to the proof of these statements. Notice that Lemma \ref{1.6}
implies Lemma \ref{1.5} that in turn implies Lemma \ref{1.4}. 

\medskip

{\it Proof of the lemmas \ref{1.1} to \ref{1.6}.}
\medskip

We choose the cohomological index of $\C^{\bullet}_{\im}D_{\bullet}$ as
row index and the homological index as column index.
The two spectral sequences coming from the horizontal and vertical
filtrations of this double complex have first terms:
$^{h}_{1}E^{i}_{j}= \C^{i}_{\im}(H_{j})$,
$^{h}_{2}E^{i}_{j}=H^{i}_{\im}(H_{j})$ and 
$^{v}_{1}E^{i}_{j}= H^{i}_{\im}(D_{j})$.

We therefore note that $(^{v}_{1}E^{i}_{j})_{\mu}=0$ for $\mu >\d_{p}$
if $i-j=p\geq 0$ and for $\mu >\e_{q}$ if $i-j=-q\leq 0$.

This implies that $H^{p}({\rm
  Tot}(\C^{\bullet}_{\im}D_{\bullet}))_{\mu}=0$ for $p\geq 0$ and $\mu 
  >\d_{p}$ and  $H^{-q}({\rm
  Tot}(\C^{\bullet}_{\im}D_{\bullet}))_{\mu}=0$ for $q\geq 0$ and $\mu
  >\e_{q}$ (with the convention ${\rm
  Tot}(C_{\bullet}^{\bullet})^{\ell}=\bigoplus_{i-j=\ell}C^{i}_{j}$).

We will now combine these vanishing with the study of the further
steps in the spectral sequences coming from the horizontal filtration
to derive the six lemmas. 

For the first three lemmas, we assume that condition $\D_{1} (\tau )$ is
satisfied and concentrate on the bottom right part of
the diagram at step 2 and relabel the second differentials that we
will use by setting  $\xi :={^{h}_{2}d^{\tau -2}_{0}}$,  $\a
:={^{h}_{2}d^{\tau -3}_{0}}$,  $\b :={^{h}_{2}d^{\tau -2}_{1}}$,  
$\g :={^{h}_{2}d^{\tau -4}_{0}}$ (see the diagram below, where the dotted arrows show direction of differentials at step 3). 

$$
\xymatrix{
\cdots&\cdots&\cdots&\cdots&^{h}_{2}E^{\tau -4}_{0}
\ar^{\g}[ddl]\\ 
\cdots&\cdots&\cdots&\cdots&^{h}_{2}E^{\tau -3}_{0}
\ar^{\a}[ddl]\ar@{..>}[dddll]\\ 
\cdots&\cdots&\cdots&^{h}_{2}E^{\tau -2}_{1}\ar^(.3){\b}[ddl]\ar@{..>}[dddll]
&^{h}_{2}E^{\tau -2}_{0}\ar^{\xi}[ddl]\ar@{..>}[dddll]\\
\cdots&\cdots&\cdots&^{h}_{2}E^{\tau -1
}_{1}\ar[ddl]\ar@{..>}[dddll]&^{h}_{2}E^{\tau -1
}_{0}\\ 
\cdots&\cdots&\, ^{h}_{2}E^{\tau }_{2}\ar[ddl]\, 
&^{h}_{2}E^{\tau}_{1}&^{h}_{2}E^{\tau}_{0}\\
\cdots&^{h}_{2}E^{\tau +1 }_{3}&^{h}_{2}E^{\tau +1 }_{2}
&0&^{h}_{2}E^{\tau +1}_{0}\\
\cdots&^{h}_{2}E^{\tau +2}_{3}&0&0&^{h}_{2}E^{\tau +2}_{0}\\
&\vdots&\vdots&\vdots&\vdots\\
}
$$

We set $\theta :={^{h}_{3}d^{\tau -3}_{0}}$.

As  for $i\geq 0$ and $\ell \geq 2$, $^{h}_{2}E^{\tau +\ell
  +i-1}_{\ell -1}=H^{\tau +\ell
  +i-1}_{\im}(H_{\ell -1})=0$, it follows that also
  $^{h}_{\ell}E^{\tau +\ell +i-1}_{\ell -1}=0$ in this case, and
  therefore ${^{h}_{\ell}d^{\tau +i-1}_{0}}:
  {^{h}_{\ell}E^{\tau +i-1}_{0}}\lra 
{^{h}_{\ell}E^{\tau +\ell +i-1}_{\ell -1}}$ is the zero map.

This implies that ${^{h}_{\ell}E^{\tau +i-1}_{0}}\simeq
{^{h}_{2}E^{\tau +i-1}_{0}}=H^{\tau +i-1}_{\im}(H_{0})$ for any
$\ell \geq 2$ and $i\geq 0$, and proves Lemma \ref{1.1}.

For Lemma \ref{1.2} first notice that $\coker (\xi)\simeq {^{h}_{3}E^{\tau
  }_{1}}\simeq {^{h}_{\infty}E^{\tau }_{1}}$.
Moreover if $\D_{2} (\tau )$ is satisfied, for $\ell \geq 2$,
  $^{h}_{2}E^{\tau +\ell -1}_{\ell}=H^{\tau +\ell -1
  }_{\im}(H_{\ell})=0$ and therefore ${^{h}_{3}E^{\tau -2}_{0}}\simeq
{^{h}_{\infty}E^{\tau -2}_{0}}$,  ${^{h}_{3}E^{\tau -1}_{1}}\simeq
{^{h}_{\infty}E^{\tau -1}_{1}}$ and ${^{h}_{4}E^{\tau }_{2}}\simeq
{^{h}_{\infty}E^{\tau }_{2}}$. The first isomorphism finishes to prove
  Lemma \ref{1.2}, and the last two isomorphisms imply Lemma \ref{1.3}
  (a). Furthermore, if  $\D_{3} (\tau )$ is satisfied then one also has   
 ${^{h}_{4}E^{\tau -3}_{0}}\simeq
 {^{h}_{\infty}E^{\tau -3}_{0}}$ and 
 ${^{h}_{3}E^{\tau -2}_{1}}\simeq
 {^{h}_{\infty}E^{\tau -2}_{1}}$ because
 $^{h}_{2}E^{\tau +\ell -1}_{\ell +1}=H^{\tau +\ell -1
  }_{\im}(H_{\ell +1})=0$ for $\ell \geq 2$, which implies Lemma \ref{1.3} (b).
\medskip

We now turn to Lemma \ref{1.6}. At step 2 we have, after relabeling the maps :
$$
\xymatrix{
\cdots&\cdots&H^{0}_{\im}(H_{3})\ar^{\phi^{3}}[ddl]&H^{0}_{\im}(H_{2})\ar^{\phi^{2}}[ddl]
\ar@{..>}^{}[dddll] 
&H^{0}_{\im}(H_{1})\ar^{\phi^{1}}[ddl]\ar@{..>}^{}[dddll]&H^{0}_{\im}(H_{0})
\ar^{\a^{0}}[ddl]\ar@{..>}^{}[dddll]\\
\cdots&\cdots&H^{1}_{\im}(H_{3})\ar_{\psi^{3}}[ddl]&H^{1}_{\im}(H_{2})\ar_{\psi^{2}}[ddl]
&H^{1}_{\im}(H_{1})\ar_{\psi^{1}}[ddl]&H^{1}_{\im}(H_{0})\ar^{\a^{1}}[ddl]\ar@{..>}^{}[dddll]\\
\cdots&H^{2}_{\im}(H_{4})&H^{2}_{\im}(H_{3})&H^{2}_{\im}(H_{2})
&H^{2}_{\im}(H_{1})\ar_(.4){\b^{0}}[ddl]&H^{2}_{\im}(H_{0})
\ar^{\a^{2}}[ddl]\ar@{..>}^{}[dddll]\\
\cdots&H^{3}_{\im}(H_{4})&H^{3}_{\im}(H_{3})&H^{3}_{\im}(H_{2})
&H^{3}_{\im}(H_{1})\ar_(.4){\b^{1}}[ddl]&H^{3}_{\im}(H_{0})
\ar^{\a^{3}}[ddl]\ar@{..>}^{}[dddll]\\
\cdots&0&0&H^{4}_{\im}(H_{2})
&H^{4}_{\im}(H_{1})\ar_(.4){\b^{2}}[ddl]&H^{4}_{\im}(H_{0})\ar^{\a^{4}}[ddl]\\
\cdots&0&0&H^{5}_{\im}(H_{2})
&H^{5}_{\im}(H_{1})&\vdots\\
\cdots&0&0&H^{6}_{\im}(H_{2})
&H^{6}_{\im}(H_{1})&\vdots\\
&\vdots&\vdots&\vdots&\vdots&\vdots\\
}
$$
\medskip
The dotted arrows shows direction of maps at step 3.

This spectral sequence abouts at step 4 because
${^{h}_{\ell}d^{i}_{j}}=0$ for any $i,j$ if $\ell \geq 4$. The
estimates we proved on the degree where the graded components of the
corresponding total complex vanishes therefore implies Lemma \ref{1.6}.\fini\medskip

These lemmas will be used frequently in the case where $D_{\bullet}$
is of the form $F_{\bullet}\otimes_{R}M$, where $F_{\bullet}$ is a
graded complex of finitely generated free $R$-modules, with $F_{i}=0$
for $i<0$ and $M$ is a graded $R$-module. Each $F_{i}$ is of the form
$\oplus_{j}R[-j]^{\b_{ij}}$ and we set $b_{i}(F_{\bullet}):=\max\{ j\
|\ \b_{ij}\not= 0\}$ if $F_{i}\not= 0$ and
$b_{i}(F_{\bullet}):=-\infty$ else. In this context, the following 
two easy lemmas will be of use:\medskip

\bl\label{1.7}  Let $F_{\bullet}$ be a graded complex of finitely
generated free $R$-modules and $M$ a
finitely generated graded $R$-module. Then  
$$
a_{p}(F_{i}\otimes_{R}M)=a_{p}(M)+b_{i}(F_{\bullet}).
$$
\el

\bl\label{1.8} Let $F_{\bullet}^{1},\ldots ,F_{\bullet}^{s}$ be
graded complexes of finitely generated free $R$-modules, with
$F_{i}^{j}=0$ for $i<0$ and any $j$, and  $F_{\bullet}$ be the tensor
product of these complexes. Then
$$
b_{\ell}(F_{\bullet})=\max_{i_{1}+\cdots +i_{s}=\ell}\{
b_{i_{1}}(F_{\bullet}^{1})+
\cdots + b_{i_{s}}(F_{\bullet}^{s})\} .
$$
\el

\section{Regularity of Frobenius powers}

Throughout this section $S$ is a standard graded ring over a field of characteristic $p>0$,
$\F$ is the Frobenius functor and $\nreg (S):=\{ \ip\in \spec (S)\ \vert\ S_{\ip}\ \hbox{is not regular}\} $.We keep the notations of sections 2 and 3.\medskip

Let $S^{[e]}$ be the ring $S$ with $S$-module structure
given by the iterated Frobenuis map $\F^{e}$. 

\bl\label{2.2}  For any  $S$-module $M$, $\tor_{i}^{S}(M,S^{[e]})$ is 
supported in $\nreg (S)\cap\supp (M)$ for any $i>0$ and any $e$.
\el

{\it Proof.} Recall that $\tor_{i}^{S}(M,S^{[e]})_{\ip}\simeq
\tor_{i}^{S_{\ip}}(M_{\ip},S_{\ip}^{[e]})$. Now $\tor_{i}^{S_{\ip}}
(M_{\ip},S_{\ip}^{[e]})=0$ for any $i$  if $\ip$ is
not in the support of $M$ and if $S_{\ip}$ is regular, 
$\F$ is exact in $S_{\ip}$ and therefore
$\tor_{i}^{S_{\ip}}(M_{\ip},S_{\ip}^{[e]}))=0$ for $i>0$. \fini\medskip  

\bt\label{2.3} Let $S$ be a standard graded ring over a field 
of characteristic $p>0$ and $M$ a finitely generated graded $S$-module. Assume 
that $\dim (\nreg (S)\cap \supp (M))\leq 1$ and set
$b_{i}^{S}(M):=\max\{ j\ |\ \tor_{i}^{S}(M,k)_{j}\not= 0\}$ if
$\tor_{i}^{S}(M,k)\not= 0$,  and $b_{i}^{S}(M):=-\infty$ else. Then, 
$$
\begin{array}{rcl}
\reg (\F^{e}M)&\leq& \max_{0\leq i\leq j\leq \dim S}\{p^{e}b_{i}^{S}(M)+a_{j}(S)+j-i\}\\
&\leq& \reg (S)+\max_{0\leq i\leq \dim S}\{ p^{e}b_{i}^{S}(M)-i\} .\\
\end{array}
$$
\et

{\it Proof.} Let $F_{\bullet}$ be a minimal graded free 
$S$-resolution of $M$, therefore $F_{i}=\oplus_{j}S[-j]^{\b_{ij}}$, where
$\b_{ij}=\dim_{k}\tor_{i}^{S}(M,k)_{j}$. 
The complex $\F^{e}F_{\bullet}$ is a complex of graded free 
$S$-modules and $\F^{e}F_{i}=\oplus_{j}S[-jp^{e}]^{\b_{ij}}$.
The hypothesis implies that $H_{i}(\F^{e}F_{\bullet})\simeq \tor_{i}^{S}(M,S^{[e]})$ 
has dimension at most $1$ for $i>0$ by Lemma \ref{2.2}. Therefore it follows from Lemma \ref{1.7} and Lemma \ref{1.1} (or
Lemma \ref{1.4})  that  
$$
a_{\ell}(\F^{e}M)\leq \max_{0\leq i\leq \dim S-\ell} \{ a_{i+\ell}(S)+p^{e}b_{i}^{S}(M)\}
, \forall \ell \geq 0,
$$
which proves our claim.\fini\medskip

\bco
Let $S$ be a standard graded ring over a field 
of characteristic $p>0$ and $M$ a finitely generated graded $S$-module. If 
$\dim (\nreg (S)\cap \supp (M))\leq 1$, and $S$ is not regular, then
$$
 \reg  (\F^{e}M)\leq \reg (S)+p^{e}\left( \reg (M)+\left\lfloor {{\dim S}\over{2}}\right\rfloor (\reg (S)-1)+\dim S\right) .
 $$
\eco

{\it Proof.} By Proposition \ref{A.9}, $b_{i}^{S}(M)\leq \reg (M)+\lfloor i/2\rfloor (\reg (S)-1)+i$\fini\medskip

This theorem gives a positive answer to question (2) in the
 introduction of \cite{Ka} when
 $\dim (\nreg (S)\cap \supp (M))\leq 1$ by providing the bound 
$$
\reg (\F^{e}M)\leq \reg (S)+p^{e}\max_{0\leq i\leq \dim S}\{b_{i}^{S}(M)\} .
$$

Notice that Lemma \ref{1.4} also shows that $\reg (\tor_{i}^{S}(M,S^{[e]}))\ll
p^{e}$ for any $i$ under the hypotheses of Theorem \ref{2.3}, and  
Lemma \ref{1.5} tells us that  $\reg (\F^{e}M)\ll p^{e}$ in
the case where $\nreg (S)\cap \supp (M)$ is of dimension two if
and only if $a_{2}(\tor_{1}^{S}(M,S^{[e]}))\ll p^{e}$ in this
case. The question of Katzman may naturally be extended by asking:
\medskip
If $F_{\bullet}$ is a graded complex of finite free $S$-modules, does one have
$\reg (H_{i}(\F^{e}F_{\bullet}))\ll p^{e}$ for any $i$ ?
\medskip
The theorem also gives an affirmative answer to Conjecture 4 in \cite{Ka} in the
special case below, where $R$ is a polynomial ring over a field of characteristic $p$,\medskip

\bco\label{2.4} Let $I,J\subseteq R$ be two homogeneous ideals and 
$L_{i}$ be the ideal generated by $X_{n},\ldots ,X_{n-i+1}$ for $i\geq 1$. Assume that 

{\rm (1)} $\proj (R/J)$ is regular outside finitely many points,\smallskip

{\rm (2)} for $i\geq 1$,  $\proj (R/(J+L_{i}))$ is regular and
$X_{n-i}$ is not a zero divisor in $R/(J+I+L_{i})^{sat}$.\smallskip

Then there exists an integer $C$ such that the Gr{\"o}bner basis of $I^{[p^{e}]}+J$ for 
the reverse lexicographic order is generated in degrees at most 
$p^{e}C$ for any $e>0$.
\eco

{\it Proof.} Condition (2) implies that the kernel of the multiplication by 
$X_{n-i}$ in $R/(J+I^{[p^{e}]}+L_{i})$ is of finite length for any $i\geq 0$ (setting $L_0:=(0)$). 
Thus $\reg (I^{[p^{e}]}+J)=\reg (in_{rev-lex}(
I^{[p^{e}]}+J))$ so that we may apply Theorem \ref{2.3} to $S:=R/J$ to get our claim.
\fini\medskip

Notice that, when condition (1) is satisfied, Bertini theorem implies
that, if $k$ is infinite, there exists $L_{i}$'s such that  $\proj
(R/(J+L_{i}))$ is regular  for $i\geq 1$. These $L_{i}$'s being
chosen, there exists a non empty open subset of the linear group such
that condition (2) is satisfied after applying a linear transformation
corresponding to a point in this subset.

\section{Regularity of  Tor modules}\medskip

Throughout this section we keep the notations of sections 2 and 3.\medskip

\bl\label{optim}
Let $S$ be a Noetherian standard graded algebra, $\im :=S_+$, $F^{\bullet}$ be a graded complex 
of finitely generated free $S$-modules, with $F^{i}=0$ for $i<0$,  $M$ be a finitely generated graded
$S$-module and $T^{\bullet}:=\C^{\bullet}_{\im}M\otimes_{S}F^{\bullet}$. 
Then,  for any $\ell$, $H^{\ell}(T^{\bullet})_{\mu}=0$ for $\mu\gg 0$. 
Consider  the two following conditions :

\smallskip
{\rm (A)} either $H^{p}_{\im}(H^{q}(F^{\bullet}\otimes_{S}M))=0$ for
$p\geq 2$ and all $q$, or there exists $q_{0}$ such that 
$H^{p}_{\im}(H^{q}(F^{\bullet}\otimes_{S}M))=0$  for $q<q_{0}$ and
$p\geq 2$, and $H^{q}(F^{\bullet}\otimes_{S}M)=0$ for $q>q_{0}$.
 
 \smallskip
{\rm (B)} $S_0$ is local with maximal ideal $\im_0$, the complex 
$F^{\bullet}\otimes_{S}S/(\im_0 +\im )$ has zero differentials and 
$B_{0}:=\sup_{q}\{ b_{0}(F^{q})+q\}<\infty$.

\smallskip
Then, 

{\rm (i)}
$$
\fin (H^{\ell}(T^{\bullet}))
\leq \max_{p+q=\ell}\{ a_{p}(H^{q}(F^{\bullet}\otimes_{S}M))\} 
$$
and equality holds if condition {\rm (A)} is satisfied.

\smallskip
{\rm (ii)}
$$
\fin (H^{\ell}(T^{\bullet}))
\leq \max_{p+q=\ell}\{ a_{p}(M)+b_{0}(F^{q})\} .
$$
If further condition {\rm (B)} is satisfied, then 
$$
\max_{j}\{ \fin (H^{j}(T^{\bullet}))+j\} 
=\reg (M)+\max_{q}\{ b_{0}(F^{q})+q\} .
$$

\smallskip
{\rm (iii)} If  conditions {\rm (A)} and  {\rm (B)} hold, then 
$$
\max_{i}\{ \reg (H^{i}(F^{\bullet}\otimes_{S}M))+i\} =\reg (M)+\max_{q}\{ b_{0}(F^{q})+q\} .
$$

\el

{\it Proof.} The homology of $T^{\bullet}$ is the aboutment of two spectral
sequences that have as first terms :
$$
{^{v}_{1}E}^{p,q}=H^{p}_{\im}(M)\otimes_{S}F^{q},
\quad\quad{^{v}_{2}E}^{p,q}=H^{q}(H^{p}_{\im}(M)\otimes_{S}F^{\bullet}),
$$
$$
{^{h}_{1}E}^{p,q}=\C^{p}_{\im}(H^{q}(M\otimes_{S}F^{\bullet})),
\quad\quad{^{h}_{2}E}^{p,q}=H^{p}_{\im}(H^{q}(M\otimes_{S}F^{\bullet})).
$$

For (i), we consider the spectral sequence obtained with the horizontal
filtration. The inequality immediately follows because
$({^{h}_{\infty}E}^{p,q})_{\mu}$ is a subquotient of
$({^{h}_{2}E}^{p,q})_{\mu}$ for any $p$, $q$ and $\mu$. Also, if  condition (A) is
satisfied, then ${^{h}_{\infty}E}^{p,q}\simeq
{^{h}_{2}E}^{p,q}$ for every $p$ and $q$, which shows the equality.\medskip

For (ii) we consider the spectral sequence obtained with the vertical
filtration. The inequality immediatly follows, because $\fin
({^{v}_{1}E}^{p,q})=a_{p}(M)+b_{0}(F^{q})$.

If condition (B) is satisfied, set $q':=\min\{ q\ \vert\ b_{0}(F^{q})+q=B_{0}\}$,
$p':=\max\{ p\ \vert\ a_{p}(M)+p=\reg (M)\}$ and $\mu
:=a_{p'}(M)+b_{0}(F^{q'})=\reg (M)+B_{0}-p'-q'$.

We will now show that $0\not= ({^{v}_{2}E}^{p',q'})_{\mu}\simeq ({^{v}_{\infty}E}^{p',q'})_{\mu}$. For
any $r\geq 0$, ${^{v}_{r+1}E}^{p',q'}$ is isomorphic to the homology of the
sequence
$$
\xymatrix{
{^{v}_{r}E}^{p'+r-1,q'-r}\ar^(.55){{^{v}_{r}d}^{p'+r-1,q'-r}}[rr]
&&{^{v}_{r}E}^{p',q'}\ar^(.45){{^{v}_{r}d}^{p',q'}}[rr]
&&{^{v}_{r}E}^{p'-r+1,q'+r}\\}
$$

Notice that $({^{v}_{1}E}^{p'-r+1,q'+r})_{\mu}\simeq
(H^{p'-r+1}_{\im}(M)\otimes_{S}F^{q'+r})_{\mu}$  is zero for $r\geq
1$ because 
$$
a_{p'-r+1}(M)+b_{0}(F^{q'+r})
\leq (\reg (M)-p'+r-1)+(B_{0}-q'-r)
\leq \mu -1.
$$
Recall that $a_{p}(M)\leq \reg (M)-p-1$ for $p>p'$ and $b_{0}(F^{q})\leq B_{0}-q-1$ for $q<q'$. This shows that $({^{v}_{1}E}^{p'+r-1,q'-r})_{\mu}=0$ for $r\geq 2$ because 
$$
a_{p'+r-1}(M)+b_{0}(F^{q'-r})
\leq (\reg (M)-p'-r)+(B_{0}-q'+r-1)
\leq \mu -1
$$

It follows that 
$$
({^{v}_{\infty}E}^{p',q'})_{\mu}\simeq ({^{v}_{2}E}^{p',q'})_{\mu}\simeq \coker (H^{p'}_{\im}(M)\otimes
F^{q'-1}\lra H^{p'}_{\im}(M)\otimes F^{q'})_{\mu}
$$

Recall that  $b_{0}(F^{q'-1})\leq b:=B_{0}-q'$. 

If $b_{0}(F^{q'-1})< b$, then $(H^{p'}_{\im}(M)\otimes
F^{q'-1})_{\mu}=0$ hence $({^{v}_{\infty}E}^{p,q})_{\mu}\simeq (H^{p'}_{\im}(M)\otimes F^{q'})_{\mu}=({^{v}_{1}E}^{p',q'})_{\mu}$. 

Else, set $F^{q'-1}=S[-b]^{\a}\oplus G^{q'-1}$ and $F^{q'}=S[-b]^{\b}\oplus G^{q'}$ with $G^{q'-1}$
and $G^{q'}$ free,  $b_{0}(G^{q'-1})<b$ and $b_{0}(G^{q'})<b$. Let $\phi$ be the map from
$S[-b]^{\a}$ to $S[-b]^{\b}$ induced by $F^{q'-1}\lra F^{q'}$, represented by a matrix $\Phi$ of elements 
in $\im_0$ as condition (B) is satisfied. 
Recall that $H:=H^{p'}_{\im}(M)_{a_{p'}(M)}$ is a finite non zero $S_0$-module, hence
$$
({^{v}_{\infty}E}^{p',q'})_{\mu}\simeq \coker (H^{\a}
\buildrel{\Phi}\over\lra H^{\b})\not= 0
$$
 by Nakayama's lemma.

Finally $\fin (H^{p'+q'}(T^{\bullet}))\geq \fin ({^{v}_{\infty}E}^{p',q'})=\mu =\reg (M)+B_{0}-p'-q'$, 
which proves our claim.\fini\medskip

Recall that for a finitely generated graded $S$-module $M$, $b_{i}^{S}(M):=\max\{ j\ \vert\ \tor_{i}^{S}(M,S_0 )_{j}\not= 0\}$ if $\tor_{i}^{S}(M,S_0 )\not= 0$ and $b_{i}^{S}(M):=-\infty$ else,
$\reg_{j}^{S}(M):=\max_{i\leq j}\{ b_{i}^{S}(M)-i\}$ and $\reg^{S}(M):=\sup_{j}\{ \reg_{j}^{S}(M)\}$. 

\bl\label{torbetti}
Let $S$ be a  standard graded ring over a Noetherian local ring, $M$
be finitely  generated graded $S$-module of finite projective dimension
and $F_{\bullet}$ be a minimal graded free $S$-resolution of $M$. Then,
$$
\reg^{S}(M)=\max_{i}\{ b_{0} (F_{i})-i \} .
$$
\el

{\it Proof.}  Let $p:=\pd M$. Applying
Lemma \ref{optim} with $F^{\bullet}:=F_{p-\bullet}$ and $M:=S_0$, so that conditions (A)  and (B) are satisfied, we obtain that
$$
\reg^{S}(M)+p=\max_{j}\{ a_{0}(\tor_{p-j}^S(M,S_0))+j\} =\max_{q}\{ b_{0}(F_{p-q})+q\}
$$
as $\reg (S_0)=0$, and the conclusion follows.\fini

\bt\label{regfpd} 
Let $S$ be a  standard graded ring over a Noetherian local ring and $M$
be finitely  generated graded $S$-module of finite projective dimension.
Then,
$$
\reg (M)=\reg^{S}(M)+\reg (S).
$$
\et

{\it Proof.} Let $F_{\bullet}$ be a minimal graded free $S$-resolution of $M$ and $p$ its length (the projective dimension of $M$).  Applying
Lemma \ref{optim} with $F^{\bullet}:=F_{p-\bullet}$ and $M:=S$, so that conditions (A) (with $q_0:=p$) and (B) are satisfied,
we obtain that 
$$
\max_{j}\{ a_{j}(M)+j+p\} =\reg (S)+\max_{q}\{ b_{0}(F_{p-q})+q\} ,
$$
and the conclusion follows from Lemma \ref{torbetti}.
\fini

\bco
Let $S$ be a  standard graded ring over a Noetherian local ring and $M$
be a finitely  generated graded $S$-module of finite projective 
dimension. Then 

(i) $\reg (M)\geq \indeg (M)+\reg (S)$,

(ii) $M$ has a linear free $S$-resolution if and only if $\reg (M)-\indeg (M)=\reg (S)$.
\eco

\bco\label{regfpd-nonlocal} 
Let $S$ be a  standard graded ring over a Noetherian ring and $M$
be finitely  generated graded $S$-module of finite projective dimension.
Then,
$$
\reg (M)\leq \reg^{S}(M)+\reg (S).
$$
\eco

{\it Proof.} Indeed, for any finitely generated graded $S$-module
$N$, $\reg (N)=\max_{\ip \in \spec (S_0)}\{ \reg (N\otimes_{S_0}(S_0)_\ip )\}$ and $b_i^{S}(N)=\max_{\ip \in \spec (S_0)}\{ b_i^{S\otimes_{S_0}(S_0)_\ip} (N\otimes_{S_0}(S_0)_\ip )\}$. Hence,
$$
\reg (M)=\max_{\ip \in \spec (S_0)}\{  \reg^{S\otimes_{S_0}(S_0)_\ip}(M\otimes_{S_0}(S_0)_\ip )+\reg (S\otimes_{S_0}(S_0)_\ip )\} \leq \reg^{S}(M)+\reg (S).
$$
\fini

\bp\label{regpolring}
Let $S$ be a  Noetherian standard graded ring over a local ring, $M$
be finitely  generated graded $S$-module and $F_\bullet$ a minimal graded free 
$S$-resolution of $M$. Then $b_{0} (F_{i})\leq \max_{j\leq i}\{ b_{j}^{S}(M)\}$ for any $i$. 
If furthermore $S_0$ has projective dimension $p<\infty$ over $S$, then   
$$
\reg (M)=\reg ^S (M)=\max_{i}\{ b_{0} (F_{i})-i \}=\max_{i\leq p}\{ b_{0} (F_{i})-i \} =\reg_p^S (M).
$$
\ep

{\it Proof.} Let $k$ be the residue field of $S_0$. The change of ring spectral sequence
$\tor_p^{S_0}(\tor_q^S(S_0,M),k)\Rightarrow \tor_{p+q}^S(M,k)$ shows that 
$b_{0} (F_{i})\leq \max_{p+q=i} \{ b_q^S(M)\} = \max_{j\leq i}\{ b_{j}^{S}(M)\}$. 

Let $F'_\bullet$  be a minimal graded free 
$S$-resolution of $S_0$ and set $b(S_0):=\max_{i}\{ b_{0} (F'_{i})-i \} $. 
By 
Theorem \ref{regfpd}, $\reg^{S}(S_0)+\reg (S)=\reg (S_0)=0$, hence
$\reg (S)=\reg^{S}(S_0)=b(S_0)=0$ by Lemma \ref{torbetti}. 

Notice that $\reg (\tor_i^S(M,S_0))=b_i^S(M)$. Therefore, applying
Lemma \ref{optim} with $F^{\bullet}:=F'_{p-\bullet}$ and $M:=M$, it follows that 
$\reg^S(M)=\reg (M)+b(S_0)=\reg (M)$. 

As $b_{j}^{S}(M)=-\infty$ for $j>p$, 
$$
\max_{i}\{ b_{0} (F_{i})-i \}\leq \max_{j\leq i}\{ b_{j}^{S}(M)-i \}=\reg_p^S(M)\leq \reg^S (M)=\reg (M).
$$
 On the other hand,
$\reg (M)\leq \reg (S)+\max_{i\leq \chdim S}\{ b_{0} (F_{i})-i \}$ by Lemma \ref{1.1} applied to $D_\bullet :=F_\bullet$.
As $\reg (S)=0$, our claim follows.\fini

\bt\label{regtor} 
Let $S$ be a standard graded ring over a Noetherian local ring $(S_0,\im_0 )$ and $M_{1},\ldots ,M_{s}$
be finitely   generated graded $S$-modules. If 

(1) $\tor_{i}^{S}(M_{1},\ldots ,M_{s})\otimes_{S_0}S_0/\im_0$ is supported in dimension at most one for any $i>0$, and 

(2) at least $s-1$ modules  among the $M_{i}$'s have finite projective dimension over $S$, then 

$$
\max_{i}\{ \reg (\tor_{i}^{S}(M_{1},\ldots ,M_{s}))-i\} = \sum_{j}\reg (M_{j})-(s-1)\reg (S).
$$
\et

{\it Proof.} We may assume that $M_i$ has finite projective dimension for $i\geq 2$. Let $F_{\bullet}$ be
the tensor product of minimal free $S$-resolutions of $M_2,\ldots ,M_s$. We set
$p:=\sum_{i=2}^{s}\pd (M_i)$ and apply Lemma \ref{optim} with $F^{\bullet}:=F_{p-\bullet}$ and $M:=M_1$.
Condition (A)  (with $q_0:=p$) and (B) are satisfied due to Proposition \ref{cohdim}. It follows that
$$
\max_{i}\{ \reg (H^{i}(F^{\bullet}\otimes_{S}M_1))+i\} = \reg (M_1)+\max_{q}\{ b_{0}(F_{p-q})+q\}.
$$
The left hand side is equal to $\max_{i}\{ \reg (\tor_{i}^{S}(M_{1},\ldots ,M_{s}))-i\}+p$ and by Lemma 
\ref{1.8} and 
 Lemma \ref{torbetti}  we can rewrite the right hand side as
$$
\reg (M_1)+\max_{i_2,\cdots ,i_s}\{ b_{i_2}^{S}(M_2)+\cdots +b_{i_s}^{S}(M_s)
-(i_2+\cdots +i_s)\} +p.
$$
This is equal to $p+\reg (M_1)+\sum_{i=2}^{s}\reg^{S}(M_i)$, and the conclusion follows from Theorem \ref{regfpd}.\fini

\bco\label{regtorex}  
Let $S$ be a standard graded ring over a Noetherian local ring $(S_0,\im_0 )$ and $M,N$
be finitely  generated graded $S$-modules.  If $M$ or $N$ have finite projective dimension  and 
$\tor_{i}^{S}(M,N)\otimes_{S_0}S_0/\im_0$ is supported in dimension at most one for $i\geq 1$,  then 
$$
\max_{i}\{ \reg (\tor_{i}^{S}(M,N))-i\} =\reg (M)+\reg (N)-\reg (S).
$$
\eco

\brm
Notice that this result doesn't hold without an hypothesis on the projective dimension of
$M$ or $N$. For instance if $M=N=S_0=S_0/\im_0$, the right hand side is negative, unless $S$ is a polynomial ring,
and the left hand side is $0$ if $S$ is Koszul and $+\infty$ else by \cite{AP}. 
\erm

\brm
Over a non local Noetherian ring $S_0$, Corollary \ref{regtorex} does not hold, as the following example shows : $S:={\bf Z}[X]$,
$M:=S\oplus S/(2)[-a]$ with $a>0$ and $N:=S/(3)$, in which case $\reg (M\otimes_S N)=\reg (N)=0$, 
$\tor_i^S(M,N)=0$ for $i>0$, $\reg (M)=a$ and $\reg (S)=0$. However, one can deduce estimates by
reduction to the local case as in Corollary \ref{regfpd-nonlocal}.
\erm

\bp\label{regtorgen} 
Let $S$ be a standard graded ring over a Noetherian local ring and $M,M_{1},\ldots ,M_{s}$
be finitely   generated graded $S$-modules. 
Set $T_{i}:=
\tor_{i}^{S}(M,M_{1},\ldots ,M_{s})$, $d:=\chdim M$,
$\tau :=\max_{i>0}\{ \chdim T_{i}\}$ and 
$$
b_{\ell}:= \max_{i_{1}+\cdots +i_{s}\leq \ell}
\{ b_{i_{1}}^{S}(M_{1})+\cdots +b_{i_{s}}^{S}(M_{s})\} .
$$
Then,\smallskip  

{\rm (i)} For $p\geq \tau -1$, 
$
a_{p}(T_{0})\leq \max_{0\leq i\leq d-p}\{
a_{p+i}(M)+b_{i} \} .
$\smallskip

{\rm (ii)}  $H^{\tau -2}_{\im}(T_{0})_{\mu}\simeq
H^{\tau}_{\im}(T_{1})_{\mu}$ for
$
\mu > \max_{0\leq i\leq d-\tau +2}\{
a_{\tau +i-2}(M)+b_{i}\} .
$\smallskip

{\rm (iii)} For 
$
\mu > \max_{0\leq i\leq d-\tau +3}\{
a_{\tau +i-3}(M)+b_{i}\} ,
$
there is an exact sequence,
$$
H^{\tau -4}_{\im}(T_{0})_{\mu}\lra H^{\tau -2}_{\im}(T_{1})_{\mu}\lra 
H^{\tau}_{\im}(T_{2})_{\mu}\lra H^{\tau -3}_{\im}(T_{0})_{\mu}
\lra H^{\tau -1}_{\im}(T_{1})_{\mu}\lra 0.
$$

{\rm (iv)} If $\tau \leq 1$, 

$\quad\quad\bullet$ $a_{p}(T_{0})\leq \max_{0\leq  i\leq d-p}\{a_{p+i}(M)+b_{i}\}$
for $p\geq 0$,\smallskip

$\quad\quad\bullet$ $a_{0}(T_{q})\leq \max_{0\leq  i\leq d}\{ a_{i}(M)+b_{q+i}\}$,
for $q\geq 0$,\smallskip

$\quad\quad\bullet$ $a_{1}(T_{q})\leq \max_{0\leq  i\leq d}\{ a_{i}(M)+b_{q+i-1}\}$,
for $q\geq 1$.\smallskip
\ep

\medskip

{\it Proof.} Let $F^{i}_{\bullet}$ be a minimal finite free
$S$-resolution of $M_{i}$ and  $F_{\bullet}$ be the tensor product
over $S$ of the complexes $F_{\bullet}^{i}$. We now apply Lemma \ref{1.7}
and lemmas \ref{1.1}, \ref{1.2}, \ref{1.3} and \ref{1.4} to $F_{\bullet}\otimes_{S}M$ to get
respectively (i), (ii), (iii) and (iv), noticing  that
$H_{q}(F_{\bullet}\otimes_{S}M)\simeq T_{q}$ by Corollary \ref{A.3} and
condition $\D_{\ell}(\tau )$ is satisfied for any $\ell\geq 1$ by
Theorem \ref{rigtor} (i).\fini\medskip

\brm
(i) If $S$ is a positively graded Noetherian algebra over a local ring $(S_0,\im_0 )$ and $N$ is a finitely generated graded $S$-module, then
$\chdim N=\dim N\otimes_{S_0}(S_0/\im_0 )$. 

(ii) By Theorem \ref{rigtor}, if $S_0$ is a field, $\tau \leq \max\{ \dim \nreg (S),\dim T_{1}\}$. 
\erm

\bco\label{3.2}  Let $k$ be a field, $\Z_{1},\ldots ,\Z_{s}$ be
closed subschemes of ${\bf P}^{n}_{k}$ and let $\Z :=\Z_{1}\cap \cdots
\cap \Z_{s}$. Assume that $\reg (\Z_{1})\geq \cdots \geq \reg
(\Z_{s})$ and there exists $\sh \subseteq \Z$ of dimension at most one
such that, locally at each point of  $\Z -\sh$, $\Z$ is a proper
intersection of Cohen-Macaulay schemes. Then, 
$$
\reg (\Z )\leq \sum_{i=1}^{\min\{ n,s\}}\reg (\Z_{i}).
$$
\eco

{\it Proof.}   Let $S$ be the polynomial ring in $n+1$ variables over $k$ and 
$I_{1},\ldots ,I_{s}$ be the defining ideals $\Z_{1},\ldots ,\Z_{s}$, 
respectively. Corollary \ref{A.10} shows that the hypothesis implies that 
$\dim \tor_{1}(S/I_{1},\ldots ,S/I_{s})\leq 2$. Therefore by
Proposition \ref{regtorgen} (with $M=S$ and $M_{i}=S/I_{i}$)
$$
\begin{array}{rl}
\reg (\Z )&=\max_{p>0}\{ a_{p}(S/(I_{1}+\cdots +I_{s}))+p\}\\
&\leq \max_{p>0}\max_{i_{1}+\cdots +i_{s}=n+1-p}\{ b_{i_{1}}(S/I_{1})+\cdots  
+b_{i_{s}}(S/I_{s})+p-n-1\}\\
&\leq \max_{p>0}\max_{i_{1}+\cdots +i_{s}=n+1-p}\{
\epsilon_{i_{1}}\reg (S/I_{1})+\cdots +\epsilon_{i_{s}}\reg (S/I_{s})\}\\
\end{array}
$$
where $\epsilon_{j}=0$ if $j\leq 0$ and  $\epsilon_{j}=1$ else (note that 
$b_{0}(S/I_{i})=0$ for any $i$). The bound follows. \fini\medskip

The following corollary refines and generalizes a theorem of Sj\"ogren in \cite{Sj}, that 
in turn extended a result of Brian\c con in \cite{Bri},
\bco\label{3.3}  Let $R$ be a polynomial ring over a field, $M$ a
finitely generated graded $R$-module  and
$f_{1},\ldots ,f_{s}$ be forms of degrees $d_{1}\geq \cdots\geq
d_{s}\geq 1$. Set $M':=M/(f_{1},\ldots ,f_{s})M$, $\d':=\dim M'$ and $\d
:=\dim M$, then for $i\geq \d '-1$
$$
a_{i}(M')\leq \max \{ a_{i}(M),a_{i+1}(M)+d_{1},\ldots
,a_{i}(M)+d_{1}+\cdots +d_{\d -i}\} .
$$
Therefore if $\dim M'=1$, 
$$
\begin{array}{rcl}
\reg (M')&\leq & \max \{ a_{0}(M),a_{1}(M)+d_{1},\ldots ,a_{\d}(M)+d_{1}+\cdots +d_{\d}\}\\
&\leq & \reg (M)+\sum_{i=1}^{\delta}(d_{i}-1),\\
\end{array}
$$
and if $\dim M'=2$, 
$$
\begin{array}{rcl}
\reg (M'/H^{0}_{\im}(M'))&\leq & \max \{ a_{1}(M),a_{2}(M)+d_{1},\ldots
,a_{\d}(M)+d_{1}+\cdots +d_{\d -1}\} +1\\
&\leq & \reg (M/H^{0}_{\im}(M))+\sum_{i=1}^{\delta -1}(d_{i}-1).\\
\end{array}
$$
\eco


{\it Proof.} Note that $\tor_{i}^{R}(M,R/(f_{1}),\ldots ,R/(f_{s}))\simeq 
H_{i}(K_{\bullet}(f_{1},\ldots ,f_{s};M))$ by Corollary \ref{A.3} and \cite[1.6.6]{BH}, 
and apply Lemma \ref{1.1}.\fini\medskip

\bex\label{3.4} Let $I$ and $J$ be two homogeneous ideals of $R:=k[
X_{1},\ldots ,X_{n}]$ with $\codim (I+J)\geq \min\{ n-1,
\codim I+\codim J\}$. Assume that $R_{\ip}/I_{\ip}$ and $R_{\ip}/J_{\ip}$ are 
Cohen-Macaulay for every homogeneous prime $\ip$ supported in $V(I+J)$ such that
$\dim R/\ip \geq 2$. Then the estimate given by Corollary \ref{3.2} gives
for any $p$ 
$$
a_{p}(R/(I+J))\leq \max_{0\leq i\leq n-p}\{ b_{i}(R/I)+b_{n-p-i}(R/J)\} -n
$$
and  Corollary \ref{3.3} gives
$$
a_{p}(R/(I+J))\leq \max_{0\leq i\leq n-p}\{ b_{i}(R/I)+a_{p+i}(R/J)\} .
$$
The first formula is more symmetric, but on the other hand in the second the
roles of $I$ and $J$ may be permuted to have another estimate, and there is
fewer terms in the second maximum as $a_{p+i}(R/J)=0$ for $i>\dim R/J-p$. 
More significatively for proving bounds, if $a_{0}(M)>a_{1}(M)>\cdots 
>a_{\dim M}(M)$ then $b_{n-p-i}(M)-n\geq a_{p+i}(M)$ for any $i$, so that in 
this case the second estimate is always sharper then the first. Note also that 
the first estimate for $p=0$ is at least $\max \{ a_{0}(R/I),a_{0}(R/J)\}$ and in
the second we can take the minimum of these two terms as the only term
concerning the $\im$-primary components of $I$ and $J$.
\eex

\bco\label{regtorsing}
Let $S$ be a standard graded ring over a field and $M,M_{1},\ldots ,M_{s}$
be finitely   generated graded $S$-modules. Assume that $\reg (S)>0$ and that
$\reg (M_i)-\indeg (M_i)\geq \reg (M_{i+1})-\indeg (M_{i+1})$ for $1\leq i\leq s-1$.
Set $T_{i}:=
\tor_{i}^{S}(M,M_{1},\ldots ,M_{s})$, $d:=\dim M$ and
$\tau :=\max_{i>0}\{ \dim T_{i}\}$. 
For $p\geq \tau -1$, 
$$
a_{p}(T_{0})\leq \max_{0\leq \ell\leq d-p}\{
a_{p+\ell }(M)+C_\ell +\ell \} ,
$$
where 
$
C_\ell :=\max_{i_1+\cdots +i_s=\ell}\{ \sum_{j}\reg_{i_j}^{S}(M_j)\}
\leq \sum_{j=1}^{s}\reg (M_j)+\lfloor \ell /2\rfloor (\reg (S)-1).
$
Therefore,
$$
\reg^{p}(T_{0})\leq
\reg^{p} (M)+\sum_{j=1}^{s}\reg (M_j)+\lfloor (d-\tau +1)/2\rfloor (\reg (S)-1).
$$
\eco

{\it Proof.} By Proposition \ref{regtorgen} (i), $a_{p}(T_{0})+p\leq \max_{0\leq i\leq d-p}\{
a_{p+i}(M)+b_{i} \}$. The number $b_i$ is bounded by Proposition \ref{A.9} (iv) :
$$
\begin{array}{rcl}
b_\ell -\ell &\leq &\max_{i_1+\cdots +i_s=\ell}\{ \sum_{j}\reg_{i_j}(M_j)+\lfloor i_j/2\rfloor (\reg (S)-1)\} \\
&\leq & \lfloor \ell /2\rfloor (\reg (S)-1)+\max_{i_1+\cdots +i_s=\ell}\{ \sum_{j}\reg_{i_j}(M_j)\} .\\
\end{array}
$$
Notice that $\reg_0 (M_{i_j})=\indeg (M_{i_j})$, hence if $\ell <s$, at least $s-\ell$ in the sum
$\sum_{j}\reg_{i_j}(M_{i_j})$ are equal to  $\reg_0 (M_{i_j})=\indeg (M_{i_j})$. This proves the first estimates, which directly implies the other one.\fini

\bt\label{regintersect}
 Let $k$ be a field, $\Z_{0}, \Z_{1},\ldots ,\Z_{s}$ be closed subschemes of $\sh$, where $\sh\subset {\bf P}^{n}_{k}$ is a closed subscheme with 
 $\reg (\sh )>0$. Set $\Z :=\Z_{0}\cap \Z_{1}\cap \cdots
\cap \Z_{s}\subset \sh$ and assume that $\reg (\Z_{1})\geq \cdots \geq \reg
(\Z_{s})$ and $\reg (\Z_{i})\geq \reg (\sh )-1$ for $i\leq \min\{ s,\dim \Z_0\}$.  

If there exists $X\subseteq \Z$ of dimension at most one
such that, locally at each point of  $\Z -X$, $\sh$ is regular, every $\Z_{i}$ is Cohen-Macaulay
and the intersection is proper in $\sh$, then, 
$$
\reg (\Z )\leq \reg (\Z_0)+\sum_{i=1}^{\min\{ s,\dim \Z_0\}}\reg (\Z_{i})+\lfloor (\dim \Z_0-1)/2\rfloor (\reg (\sh )-1).
$$
\et

{\it Proof.} Let $S$ be the standard graded algebra with positive depth such that $\sh =\proj (S)$,
$A_i$ be the graded quotient of $S$ with positive depth defining $\Z_i$,  $d:=\dim A_0=\dim (\Z_0)+1$ and $e :=\min\{ d-1,s\}$. By Theorem \ref{rigtor}, $\tau \leq 2$. Therefore Corollary \ref{regtorsing} implies that
$$
\reg^{1}(B)\leq \reg (A_0)+\max_{i_1+\cdots +i_s=d-1}\{ \sum_{j}\reg_{i_j}^{S}(A_j)\}
$$
and by Example \ref{bettideaux}
$$
\reg_{i_j}^{S}(A_j)\leq  \max \{ \reg (A_j),\reg (S)-1\} +\lfloor (i_j-1)/2\rfloor ( \reg(S)-1).
$$
As $\reg_{0}^{S}(A_j)=0$ for any $j$,  $\reg (A_{1})\geq \cdots 
\geq \reg (A_{s})$ and $\reg (A_j)\geq \reg (S)-1$ for $j\leq e$, we get :
$$
\reg (\Z )=\reg^{1}(B)\leq \reg (A_0)+\sum_{i=1}^{e}\reg (A_j)+\lfloor (d-2)/2\rfloor ( \reg(S)-1)
$$
because $\max_{i_1+\cdots +i_e=d-1}\{ \sum_{j} \lfloor (i_j-1)/2\rfloor \} =\lfloor (d-2)/2\rfloor $.
\fini

\brm
In the case $\dim \Z_0< 2s$, the theorem may be refined slightly as follows :
$$
\reg (\Z )\leq \reg (\Z_0)+\sum_{i=1}^{\dim \Z_0}\reg (\Z_{i}),\quad {\rm if}\ \dim \Z_0 \leq s
$$
$$
\reg (\Z )\leq \reg (\Z_0)+\sum_{i=1}^{\dim \Z_0}\reg (\Z_{i})+(\dim \Z_0 -s)(\reg (S)-1),\quad {\rm if}\ s\leq \dim \Z_0 \leq 2s.
$$
\erm

\bex
Let $k$ be a field, $\Z_{0}, \Z_{1},\ldots ,\Z_{s}$ be closed subschemes of $\sh\subset {\bf P}^{n}_{k}$ with 
 $\reg (\sh )=1$ and let $\Z :=\Z_{0}\cap \Z_{1}\cap \cdots
\cap \Z_{s}\subset \sh$. Assume that $\reg (\Z_{1})\geq \cdots \geq \reg
(\Z_{s})$.  

If there exists $X\subseteq \Z$ of dimension at most one
such that, locally at each point of  $\Z -X$, $\sh$ is regular, every $\Z_{i}$ is Cohen-Macaulay
and the intersection is proper in $\sh$, then, 
$$
\reg (\Z )\leq \reg (\Z_0)+\sum_{i=1}^{\min\{ s,\dim \Z_0\}}\reg (\Z_{i}).
$$
\eex

\bex
Let $\Z_{1},\Z_{2}$ be two reduced  divisors of $\sh\subset {\bf P}^{n}_{k}$, where $\sh$ is irreducible of dimension $4$ with
a singular locus of dimension at most 1. If $\dim (\Z_{1}\cap \Z_{2})=2$, then :
$$
\reg (\Z_{1}\cap \Z_{2})\leq   \reg (\Z_{1})+\max\{ \reg (\Z_{2}),\reg (\sh )-1\} +\reg (\sh )-1.
$$
\eex

\section{On the regularity of ordinary powers}

In this section $S$ is a Noetherian standard graded algebra, and $R$ is 
a standard graded polynomial ring over a Noetherian ring, and we keep the notations 
of sections 2 and 3.\medskip

Recall that for a finitely generated graded $S$-module, the cohomological dimension of $M$ (with respect to $\im =S_+$) is
$\chdim M:= \max\{ i\ \vert\ H^i_{\im}(M)\not= 0\} =\max_{\ip\in \spec (S_0)}\{ \dim (M\otimes_{S_0} k(\ip ))\}$ (where $k(\ip )$ is the residue field of $(S_0)_\ip$). 

We will  control
the regularity of $I^{m}/I^{m+1}=\tor_{1}^{S}(S/I,S/I^{m})$, for
$m\geq 1$ in terms of invariants attached to $I$. We start with a result that generalizes and refines previous results in \cite{Chan}, \cite{GGP} and \cite{EHU}.\medskip

\bt\label{regpow1}
Let $I$ be an homogeneous $S$-ideal  such that $\chdim (S/I)\leq 1$. Then, for any $m\geq 0$,

(i) $a_{1}(S/I^{m+1})\leq a_{1}(S/I)+mb_{0}^{S}(I)$,

(ii) $a_{0}(S/I^{m+1})\leq \max\{ a_{0}(S/I)+b_{0}^{S}(I),a_{1}(S/I)+b_{1}^{S}(I)\} +(m-1)b_{0}^{S}(I).$

Hence, 
$$
\reg (S/I^{m+1})\leq \max\{ a_{0}(S/I)+b_{0}^{S}(I),a_{1}(S/I)+1+\reg_{1}^{S}(I)\} +(m-1)b_{0}^{S}(I). 
$$
\et

{\it Proof.} First notice that we may assume that $S_0$ is local. Furthermore 
$\chdim (\tor_{i}^{S}(S/I,S/I^{m}))\leq 1$ for all $i$ and all $m$ by \cite[2.2]{DNT} as
$\supp (\tor_{i}^{S}(S/I,S/I^{m}))\subseteq \supp (S/I)$. 

Applying Proposition \ref{regtorgen} (iv) to $M:=S/I^m$ and $M_1:=S/I$, it follows that 
$$
(*)\quad a_{0}(I^m/I^{m+1})\leq \max\{ a_{0}(S/I^m)+b_{0}^{S}(I),a_{1}(S/I^m)+b_{1}^{S}(I)\}
$$
and applying Proposition \ref{regtorgen} (iv) to $M/H^0_\im (M)=S/(I^m)^{sat}$ and $M_1$, it follows that 
$$
a_{1}((I^m)^{sat} /I(I^{m})^{sat} )\leq a_{1}(S/(I^m)^{sat})+b_{0}^{S}(I)
$$
but $H^1_\im (S /(I^{m})^{sat} )\simeq H^1_\im (S /I^{m})$ and the exact sequence 
$$
H^1_\im ((I^m)^{sat} /I(I^{m})^{sat} )\ra 
H^1_\im (S /I(I^{m})^{sat} )\ra H^1_\im (S /(I^{m})^{sat} )
$$
 together with the isomorphim $H^1_\im (S /I(I^{m})^{sat} )\simeq H^1_\im (S /I^{m+1} )$ implies that
 $$
(**)\quad  a_{1}(S/I^{m+1})\leq a_{1}(S/I^m)+b_{0}^{S}(I).
 $$
The inequalities $(*)$ and $(**)$ show (i) and (ii) by recursion on $m$.\fini

 \bco
 Under the hypotheses of Theorem \ref{regpow1}, if $S=R/J$, where $R$ is  a polynomial ring over
 $S_0$,
 $$
\reg (S/I^{m+1})
\leq \max\{ \reg (S/I),b_0^R(J)-2\} +\max \{  \reg (S/I^{sat})+1,b_{0}^{S}(I)\} +(m-1)b_{0}^{S}(I) .
$$
\eco

{\it Proof.} Recall that  $\reg_1^S(I)\leq \max\{ \reg_1^R(I+J),\reg_1^R(S)\} =\max\{ \reg_1^R(I+J),b_0^R(J)-1\}$, so that
by Theorem \ref{regpow1} and Proposition \ref{regpolring} we have,
  $$
\begin{array}{rl}
\reg (S/I^{m+1})&\leq \max\{ a_{0}(S/I)+b_{0}^{S}(I),a_{1}(S/I)+1+\reg_{1}^{S}(I)\} +(m-1)b_{0}^{S}(I) \\
 &\leq \max\{ a_{0}(S/I)+b_{0}^{S}(I),a_{1}(S/I)+1+\reg_{2}^{R}(S/I),a_{1}(S/I)+b_0^R(J)\} +(m-1)b_{0}^{S}(I)\\
\end{array}
$$
and $a_{1}(S/I)+1=\reg (S/I^{sat})$. This implies our claim.
\fini

\medskip
The next lemma will be useful to study the cases of dimension two and
three.\medskip 

\bl\label{4.2}  Let $I$ be an homogeneous  $S$-ideal 
generated in degrees $d_{1}\geq \cdots \geq d_{s}$ and $d:=\dim
(S/I)$. Assume that $I_{\ip}\subset S_{\ip}$ is a complete
intersection  for every prime $\ip\subseteq \supp (S/I)$ of maximal
dimension. Then, 
$$
a_{d}(\tor_{m}^{S}(S/I,S/I))\leq a_{d}(\bigwedge^{m}(I/I^{2}))\leq 
a_{d}(S/I)+d_{1}+\cdots +d_{m}.
$$
\el

{\it Proof.}  Set $B:=S/I$. The alternating algebra structure on
$\tor_{\bullet}^{S}(B,B)$ and the  
identification $ I/I^{2}\buildrel{\psi}\over{\lra}\tor_{1}^{S}(B,B)$ gives rise 
to a $B$-algebra homomorphism:
$$
\bigwedge \psi :\bigwedge (I/I^{2})\lra \tor_{\bullet}^{S}(B,B)
$$
which is an isomorphism  locally at every prime $\ip$ such that $\dim S/\ip =d$
by \cite[2.3.9]{BH}. Therefore  $a_{d}(\tor_{m}^{S}(B,B))\leq
a_{d}(\bigwedge^{m}(I/I^{2}))$. Also if $I=(f_{1},\ldots ,f_{s})$ with
$\deg f_{i}=d_{i}$, there is a natural onto map 
$$
T_{m}:=\bigotimes_{i=1}^{m}(f_{i},\ldots ,f_{s})/I(f_{i},\ldots ,f_{s})\lra
\bigwedge^{m}(I/I^{2})
$$
of graded modules supported in dimension at most $d$, which shows that
$a_{d}(\bigwedge^{m}(I/I^{2}))\leq a_{d}(T_{m})$. By Proposition 
\ref{regtorgen}
(i) we have $a_{d}(T_{m})\leq a_{d}(S/I)+d_{1}+\cdots +d_{m}$, because
$(f_{i},\ldots ,f_{s})/I(f_{i},\ldots ,f_{s})$ is generated in degree
at most $d_{i}$. 
This proves our result. 
\fini\medskip    

\bt\label{4.3} Let $I$ be an homogeneous $S$-ideal 
such that $\dim (S/I)=2$. Assume that $I_{\ip}\subset S_{\ip}$ is a complete 
intersection  for every prime $\ip$ such that $\dim S/\ip =2$. Set $a_{i}:=a_{i}
(S/I)$, $b_{i}:=b_{i}^{S}(I)$ and let $d_{1}\geq \cdots \geq d_{s}$ be the degrees of a 
minimal system of generators of $I$. Then, for any $m>0$
$$
\begin{array}{rl}
a_{0}(\tor_{m}^{S}(S/I,S/I))&\leq \max\{ a_{0}+b_{m-1},a_{1}+b_{m},a_{2}+b_{m+1},
a_{2}(\tor_{m+1}^{S}(S/I,S/I))\} \\
&\leq \max\{ a_{0}+b_{m-1},a_{1}+b_{m},a_{2}+b_{m+1},a_{2}+d_{1}+\cdots +d_{m+1}\} ,
\\
\end{array}
$$
$$
a_{1}(\tor_{m}^{S}(S/I,S/I))\leq \max\{ a_{1}+b_{m-1},a_{2}+b_{m}\}
$$
and
$$
a_{2}(\tor_{m}^{R}(S/I,S/I))\leq a_{2}+d_{1}+\cdots +d_{m}.
$$
\et

Note that in the case $S=R$ is a polynomial ring over a field, $\tor_{m}^{R}(R/I,R/I)=0$ and 
$b_{m}=0$ if $m>\dim R =:n$.  In this case, $a_{0}(\tor_{n}^{R}(R/I,R/I))\leq
a_{0}+b_{n-1}$.\medskip

{\it Proof.} By Lemma \ref{4.2}, $a_{2}(\tor_{m}^{R}(B,B))\leq
a_{2}+d_{1}+\cdots +d_{m}$ and the result follows from Lemma \ref{1.5} (b).\fini\medskip

From this point on, we assume that $R_0=k$ is a field.\medskip

\bt\label{rpd2}  Let $S$ be a standard graded algebra over a field, $I$ be an homogeneous $S$-ideal 
such that $\dim (S/I)=2$ and $d_{1}\geq \cdots \geq d_{s}$ be the degrees of a 
minimal system of generators of $I$. For $\ell\geq 2$, set $E_{\ell }:=\max_{i\geq j}\{ a_{i}(S/I^{\ell -1 })+b_{j}^{S}(I)\}$ and
let $T_{\ell }:=\ker (I\otimes_{S}I^{\ell -1}\ra I^{\ell })$. Then,\smallskip

(i) If $S$ is regular $E_{\ell }\leq \reg (S/I^{\ell -1})+\reg (I)$, and if $S=R/J$, where $R$ is a polynomial ring and $J\not= 0$, one has
$$
E_{\ell }
\leq  \max\{ \reg (S/I^{\ell -1})+ \reg (I),a_{2}(S/I^{\ell -1 })+b_{0}^{R}(J)+b_{0}^{R}(I)\} .
$$

(ii) $\reg (S/I^{\ell })\leq \max\{ E_{\ell},a_{2}(T_{\ell })\}$. Moreover, $\reg (R/I^{\ell })> E_{\ell}$ if and only if
$a_{2}(T_{\ell })>E_{\ell}$, in which case $\reg (R/I^{\ell })=a_{2}(T_{\ell })$.\smallskip

Assume further that $I_{\ip}\subset S_{\ip}$ is a complete 
intersection  for every prime $\ip$ such that $\dim S/\ip =2$. Then,\smallskip

(iii) $a_{2}(T_{2})\leq a_{2}(S/I)+d_{1}+d_{2}$ and $a_{2}(T_{\ell })\leq a_{2}(T_{\ell -1})+b_{0}(I)$ for $\ell >2$,\smallskip

(iv) Set $E'_3:=\max\{ E_{2}+b_{0}^{S}(I),a_{2}(S/I)+2b_{1}^{S}(I)\}$. Then,
$E'_3\leq \reg (S/I)+2\reg (I)$ if $S$ is regular. If $S=R/J$, where $R$ is a polynomial ring and $J\not= 0$, one has
$$
E'_3
\leq  \max\{ \reg (S/I)+2\reg (I),a_{2}(S/I)+b_{0}^{R}(J)+2b_{0}^{R}(I)\} .
$$
Furthermore, $\reg (S/I^{2})\leq \max\{ E_{2},a_{2}(S/I)+d_{1}+d_{2}\}$ and
$$
\reg (S/I^{\ell })\leq \max\{ E'_3,a_{2}(S/I)+2d_{1}+d_{2}\}+(\ell -3)d_{1},\quad \forall \ell \geq 3.
$$

\et

{\it Proof.} Part (i) follows from the definition of regularity in terms of graded Betti numbers over a regular ring
together with Example \ref{bettideaux}.

We now prove (ii). Let $F_{\bullet}$ be a minimal free $S$-resolution of $S/I$ and consider $D_{\bullet}^{\ell }:=F_{\bullet}\otimes_{S}
S/I^{\ell }$. One has, $H_{i}(D_{\bullet}^{\ell })=\tor_{i}^{S}(S/I,S/I^{\ell })$, in particular  $H_{0}(D_{\bullet}^{\ell })=S/I^{\ell }$,
 $H_{1}(D_{\bullet}^{\ell })\simeq I^{\ell -1}/I^{\ell }$ and $H_{2}(D_{\bullet}^{\ell })\simeq T_{\ell }$. 
 
 Set $a_{i}^{\ell }:=a_{i}(R/I^{\ell })$,  $b_{i}^{S}:=b_{i}^{S}(S/I)$,  $B_{i}^{S}:=\max_{j\leq i}\{ b_{j}^{S}\}$ and $b_{i}^{R}:=b_{i}^{R}(S/I)$.
 
 It is easily seen that $a_{2}^{\ell }\leq a_{2}^{\ell -1}+b_{0}$. 
 By Proposition \ref{regtorgen} (ii) applied to $M:=S/I^{\ell -1}$ and $M_1:=S/I$, one has $a_{1}(I^{\ell -1}/I^{\ell })\leq \max\{ a_{1}^{\ell -1}+b_{0},a_{2}^{\ell -1}+b_{1}\}$ 
 which implies that $a_{1}^{\ell }\leq \max\{ a_{1}^{\ell -1}+B_{0},a_{2}^{\ell -1}+B_{1}\}$, and 
there is a natural map 
$$
\phi^{\ell }_{\mu} : H^{0}_{\im}(I^{\ell -1}/I^{\ell })_{\mu}\lra H^{2}_{\im}(T_{\ell })_{\mu}
$$ 
which is an isomorphism for $\mu >\max\{ a_{0}^{\ell -1}+B_{0},a_{1}^{\ell -1}+B_{1},a_{2}^{\ell -1}+B_{2}\}$. In particular, it
shows that for $\mu > B_{\ell -1}$ there is an isomorphism $H^{0}_{\im}(S/I^{\ell })_{\mu}\simeq H^{2}_{\im}(T_{\ell })$.

For (iii) recall that $a_{2}(T_{2})\leq a_{2}(S/I)+d_{1}+d_{2}$ by Lemma \ref{4.2}. Notice that for $\ell \geq 2$, 
there is commutative diagram
$$
\xymatrix{
&T_{2}\otimes I^{\ell -1}\ar[r]\ar@{-->}^{\psi_{\ell }}[d]&I\otimes I\otimes I^{\ell -1}\ar[r]\ar[d]&I^{2}\otimes I^{\ell -1}
\ar[r]\ar[d]&0\\
0\ar[r]&T_{\ell }\ar[r]&I\otimes I^{\ell }\ar[r]\ar[d]&I^{\ell +1}\ar[r]\ar[d]&0.\\
&&0&0& \\}
$$
We first show that $\psi_{\ell }$ is generically onto. This is the case if the induced map 
$\tau_{\ell }: \bigwedge^{2}(I/I^{2})\otimes (I^{\ell -1}/I^{\ell })\lra T_{\ell }$ is generically onto. When 
$I$ is a complete intersection,  $I/I^{2}$ is free, $\sym^{\ell}_{S}(I/I^{2})=I^{\ell}/I^{\ell +1}$, and there is a commutative
diagram with  exact  rows
$$
\xymatrix@C=10pt{
&\bigwedge^{2}(I/I^{2})\otimes \sym^{\ell -1}_{S}(I/I^{2})\ar^{\tau_{\ell }}[d]\ar[r]&
(I/I^{2})\otimes \sym^{\ell }_{S}(I/I^{2})\ar[r]\ar^{\simeq}[d]&\sym^{\ell +1}_{S}(I/I^{2})\ar^{\simeq}[d]\ar[r] &0\\
0\ar[r]&T_{\ell }\ar[r]&(I/I^{2})\otimes (I^{\ell }/I^{\ell +1})\ar[r]&I^{\ell +1}/I^{\ell +2}\ar[r]&0\\
}
$$
showing the surjectivity of  $\tau_{\ell }$. Therefore  $\psi_{\ell}$  is indeed generically onto.

Now, as $\psi_{\ell }$ has a cokernel of dimension at most 1, 
$$
a_{2}(T_{\ell })\leq a_{2}(T_{2})+(\ell -1)b_{0}(I)\leq a_{2}+\ell d_{1}+d_{2},
$$
by Theorem 4.3 and Proposition 
\ref{regtorgen} (i). 

By (i), (ii) and (iii), for $\ell \geq 2$, $a_{2}^{\ell }\leq a_{2}^{\ell -1}+B_{0}$,
$a_{1}^{\ell }\leq \max\{ a_{1}^{\ell -1}+B_{0},a_{2}^{\ell -1}+B_{1}\}$ and therefore
$$
a_{1}^{\ell }\leq  \max\{  a_{1}+(\ell -1)B_{0},a_{2}+B_{1}+(\ell -2)B_{0}\},
$$

which shows that $a_{0}^{2}\leq \max\{ E_{2},a_{2}+d_{1}+d_{2}\}$ using (ii), and gives for $\ell >2$,
$$
\begin{array}{rl}
a_{0}^{\ell }&\leq \max\{ a_{0}^{\ell -1}+B_{0},a_{1}^{\ell -1}+B_{1},a_{2}^{\ell -1}+B_{2},a_2  +(\ell +1)B_0\}\\
&\leq \max\{ a_{0}^{\ell -1}+B_{0},\max\{ a_{1}+B_0+B_{1},a_{2}+2B_{1},a_{2}+B_0+B_{2},a_2  +3 B_0\} +(\ell -3)B_0\} \\
&\leq \max\{ a_{0}+2B_{0},a_{1}+B_0+B_{1},a_{2}+2B_{1},a_{2}+B_0+B_{2},a_2  +3 B_0\} +(\ell -3)B_0\\
\end{array}
$$
which proves our claim together with the above estimates on $a_{1}^{\ell }$ and $a_{2}^{\ell }$.
\fini\medskip

\bco\label{4.4}  Let $I$ be an homogeneous ideal of
  $R$  
such that $\dim (R/I)=2$ and $d_{1}\geq \cdots \geq d_{s}$ be the degrees of a 
minimal system of generators of $I$. Set $a_{i}:=a_{i}(R/I)$, $b_{i}:=b_{i}(I)$ and
$b_{i}':=\max_{j=1,\ldots ,i}\{ b_{j}\} $. 
Assume that $I_{\ip}\subset R_{\ip}$ is a complete 
intersection  for every prime $\ip$ such that $\dim R/\ip =2$. 
$$
\begin{array}{rl}
\reg (I^{2})&\leq \max\{ a_{0}+b_{0},a_{1}+b'_{1},a_{2}+b'_{2},a_{2}+d_{1}+d_{2}\} +1\\
&\leq \max \{ \reg (I)+\max\{ b_{0},b_{1}-1,b_{2}-2\} ,a_{2}+2b_{0}+1\}.\\
&\leq \max \{ 2\reg (I),\reg (I^{sat})+2b_{0}-2\}\\
\end{array}
$$
For $j\geq 2$, set $B:=\max \{ b'_2+b_0,2b'_1,2d_1+d_2\}$, then
$$
\begin{array}{rl}
\reg (I^{j})&\leq \max\{ a_{0}+2b_{0},a_{1}+b'_{1}+b_{0},a_{2}+B\}+(j-3)b_{0}+1\\
&\leq \max \{ 3\reg (I)+(j-3)b_{0},a_{2}+jb_{0}+1\}\\
&\leq  \max \{ 3\reg (I)+(j-3)b_{0},\reg (I^{sat})+jb_{0}-2\}\\
\end{array}
$$
\eco

\bt\label{4.5} Let $I$ be an homogeneous ideal of $R$ 
such that $\dim (R/I)=3$. Assume that $I_{\ip}\subset R_{\ip}$ is a complete 
intersection for every prime $\ip$ such that $\dim R/\ip =3$. Set $a_{i}:=a_{i}(R/I)$, 
$b_{i}:=b_{i}(I)$ and let $d_{1}\geq \cdots \geq d_{s}$ be the degrees of a 
minimal system of generators of $I$. Then,\smallskip

{\rm (i)} there exists a natural map $\psi : H^{0}_{\im}(I/I^{2})\lra 
H^{2}_{\im}(\tor_{2}(R/I,R/I))$ such that $\psi_{\mu}$ is onto for $\mu >  
\max \{ a_{1}+b_{0},a_{2}+b_{1},a_{3}+b_{2}\}$
and into for $\mu > \max \{ a_{0}+b_{0},a_{1}+b_{1},a_{2}+b_{2},a_{3}+b_{3},
a_{3}+d_{1}+d_{2}+d_{3}\}$,\smallskip

{\rm (ii)}  $H^{1}_{\im}(I/I^{2})_{\mu}=0$ for 
$\mu >  \max \{ a_{1}+b_{0},a_{2}+b_{1},a_{3}+b_{2},a_{3}+d_{1}+d_{2}\}$,\smallskip

{\rm (iii)}  $H^{2}_{\im}(I/I^{2})_{\mu}=0$ 
 for $\mu >  \max \{ a_{2}+b_{0},a_{3}+b_{1}\}$,\smallskip

{\rm (iv)} $H^{3}_{\im}(I/I^{2})_{\mu}=0$ 
 for $\mu >  a_{3}+b_{0}$.
\et

{\it Proof.} Let $F_{\bullet}$ be a minimal free resolution of $R/I$. 
As $I\subset R$ is a generically a complete intersection,
$a_{3}(\tor_{m}^{R}(R/I,R/I)\leq a_{3}+d_{1}+\cdots +d_{m}$ by Lemma
\ref{4.2}. The result then follows from Lemma \ref{1.6} applied to $D_\bullet :=F_{\bullet}\otimes_R R/I$.\fini\medskip  

As a corollary, we present an application for bounding the regularity of the module
of K\"ahler differentials associated to a standard graded algebra defining a projective surface. Similar results for the four highest
cohomology modules of the module of K\"ahler differentials can be obtained in general.

\bco\label{4.6} Let $I$ be an homogeneous ideal of 
$R$, $B:=R/I$ and $\sh :=\proj (B)$ be the corresponding 
projective scheme. Assume that $\sh$ is of dimension two and generically a complete 
intersection. Set $a_{i}:=a_{i}(B)$, $b_{i}:=b_{i}(I)$ and let 
$d_{1}\geq \cdots \geq d_{s}$ be the degrees of a 
minimal system of generators of $I$. Then if $\Om_{B}$ is the module
of K{\"a}hler differentials of $B$ (over $k$), one has\smallskip

{\rm (i)} $a_{3}(\Om_{B})\leq a_{3}+1$,  $a_{2}(\Om_{B})\leq \max \{ a_{2}+1,
a_{3}+b_{0}\}$,\smallskip

{\rm (ii)} $a_{1}(\Om_{B})\leq \max \{ a_{1}+1,a_{2}+b_{0},a_{3}+b_{1}\}$ if $\sh$
is generically reduced,\smallskip

{\rm (iii)} $a_{0}(\Om_{B})\leq \max \{ a_{0}+1,a_{1}+b_{0},a_{2}+b_{1},a_{3}+b_{2},
a_{3}+d_{1}+d_{2}\}$ if $\sh$ is a reduced complete intersection outside finitely many points.
\eco

{\it Proof.} We consider the exact sequence presenting $\Om_{B}$ and
defining $K$ :
$$
0\lra K\lra I/I^{2}\buildrel{\psi}\over{\lra}B[-1]^{n}\lra \Om_{B}
\lra 0,
$$
where $\psi$ is given by the jacobian matrix on a minimal system of generators of 
$I$. Taking cohomology on the sequence  
we get (i) from Theorem \ref{4.5} (iv). If $\sh$ is generically reduced, then
$K$ is supported in dimension at most two, and therefore (ii) follows from
Theorem \ref{4.5} (iii). In the situation of (iii) $\dim K\leq 1$, and
Theorem \ref{4.5} (ii) proves the claim.\fini

\vspace{1in}

\end{document}